\documentclass[10pt,twoside]{article}
\usepackage{fix-cm}

\usepackage[T1]{fontenc}
\usepackage[british]{babel}
\usepackage{csquotes}
\usepackage{microtype}
\usepackage[
  a4paper,
  inner=32mm,
  outer=32mm,
  top=32mm,
  bottom=32mm,
  heightrounded
]{geometry}
\usepackage{setspace}
\usepackage{amsmath,amsthm,mathtools}
\usepackage[cmintegrals,cmbraces]{newtxmath}
\usepackage{ebgaramond-maths}
\usepackage{bussproofs-extra}

\DeclareSymbolFont{cmletters}{OML}{cmm}{m}{it}
\makeatletter
\re@DeclareMathSymbol{\triangleright}{\mathbin}{cmletters}{"2E}
\re@DeclareMathSymbol{\frown}{\mathrel}{cmletters}{"5F}
\re@DeclareMathAccent{\vec}{\mathord}{symbols}{174}
\makeatother

\usepackage{enumitem}
\usepackage{float}
\usepackage{titlesec}
\usepackage{xcolor}
\usepackage{xspace}

\newcommand{\paperkeywords}{%
  reflection calculus; quantified modal logic; cut elimination;
  uniform interpolation; provability logic; conservativity results}
\newcommand{\papermsc}{%
  Primary 03F03, 03B45, 03F45; Secondary 03C40}

\usepackage[hidelinks]{hyperref}
\hypersetup{
  pdftitle={Modular proof analysis and cut elimination for RC1 and QRC1},
  pdfauthor={Cosimo Perini Brogi and Joost J. Joosten},
  pdfsubject={Cut-free proof theory for reflection calculi},
  pdfkeywords={\paperkeywords}
}
\usepackage[nameinlink,noabbrev]{cleveref}

\usepackage[
  style=numeric-comp,
  sorting=nyt,
  maxbibnames=99
]{biblatex}
\AtBeginBibliography{\footnotesize\sloppy}
\titleformat{\section}[block]
  {\normalfont\large\scshape\filcenter}{\thesection}{0.75em}{}
\titleformat{\subsection}[block]
  {\normalfont\itshape\filcenter}{\thesubsection}{0.75em}{}
\titleformat{\paragraph}[runin]
  {\normalfont\scshape}{}{0pt}{}[.]
\titlespacing*{\section}
  {0pt}{2.4ex plus .7ex minus .2ex}{1.0ex}
\titlespacing*{\subsection}
  {0pt}{1.8ex plus .5ex minus .2ex}{0.7ex}
\titlespacing*{\paragraph}
  {0pt}{1.4ex plus .4ex minus .2ex}{0.8em}

\newtheoremstyle{tufteplain}
  {0.9\baselineskip}{0.55\baselineskip}
  {\itshape}{}{\scshape}{.}{0.5em}{}
\newtheoremstyle{tuftedefinition}
  {0.9\baselineskip}{0.55\baselineskip}
  {\normalfont}{}{\scshape}{.}{0.5em}{}
\newtheoremstyle{tufteremark}
  {0.75\baselineskip}{0.45\baselineskip}
  {\normalfont}{}{\itshape}{.}{0.5em}{}

\theoremstyle{tufteplain}
\newtheorem{theorem}{Theorem}[section]
\newtheorem{proposition}[theorem]{Proposition}
\newtheorem{lemma}[theorem]{Lemma}
\newtheorem{corollary}[theorem]{Corollary}
\theoremstyle{tuftedefinition}
\newtheorem{definition}[theorem]{Definition}

\theoremstyle{tufteremark}

\crefname{theorem}{theorem}{theorems}
\Crefname{theorem}{Theorem}{Theorems}
\crefname{proposition}{proposition}{propositions}
\Crefname{proposition}{Proposition}{Propositions}
\crefname{lemma}{lemma}{lemmas}
\Crefname{lemma}{Lemma}{Lemmas}
\crefname{corollary}{corollary}{corollaries}
\Crefname{corollary}{Corollary}{Corollaries}
\crefname{definition}{definition}{definitions}
\Crefname{definition}{Definition}{Definitions}
\crefname{example}{example}{examples}
\Crefname{example}{Example}{Examples}
\crefname{remark}{remark}{remarks}
\Crefname{remark}{Remark}{Remarks}

\newcommand{\gl}{\ensuremath{\mathsf{GL}}\xspace}
\newcommand{\SGL}{\ensuremath{\mathsf{SGL}}\xspace}
\newcommand{\SGLDiamond}{\ensuremath{\mathsf{SGL}^{\Diamond}}\xspace}
\newcommand{\RC}{\ensuremath{\mathsf{RC}}\xspace}
\newcommand{\QRC}{\ensuremath{\mathsf{QRC}_1}\xspace}
\newcommand{\lang}{\ensuremath{\mathcal L^+}\xspace}
\newcommand{\qlang}{%
  \ensuremath{\mathcal{L}^{+}_{\forall}}\xspace}

\newcommand{\calc}{\ensuremath{\mathrm{GRC1}}\xspace}
\newcommand{\qcalc}{\ensuremath{\mathrm{GQRC1}}\xspace}
\newcommand{\calcFourNecStar}{%
  \ensuremath{\mathrm{GK4}^{+}_{*}}\xspace}
\newcommand{\qcalcFourNecStar}{%
  \ensuremath{\mathrm{GQK4}^{+}_{*}}\xspace}
\newcommand{\calcFourNec}{%
  \ensuremath{\mathrm{GK4}^{+}}\xspace}
\newcommand{\qcalcFourNec}{%
  \ensuremath{\mathrm{GQK4}^{+}}\xspace}

\newcommand{\subst}[2]{[#1\leftarrow #2]}
\newcommand{\FV}{\operatorname{FV}}
\newcommand{\BV}{\operatorname{BV}}
\newcommand{\Rel}{\operatorname{Rel}}
\newcommand{\Const}{\operatorname{Const}}
\newcommand{\Tm}{\operatorname{Tm}}
\newcommand{\Sub}{\operatorname{Sub}}

\newcommand{\mdepth}{\mathrm{d}_{\Diamond}}
\newcommand{\qdepth}{\mathrm{d}_{\forall}}

\setlist[itemize]{
  leftmargin=*,
  itemsep=0.35em,
  topsep=0.45em
}
\setlist[enumerate]{
  leftmargin=*,
  itemsep=0.35em,
  topsep=0.45em
}

\title{Modular proof analysis and cut elimination\\[0.1em]
  for the Reflection Calculi $\mathsf{RC}_1$ and $\mathsf{QRC}_1$}
\author{%
  Cosimo Perini Brogi\\[0.35em]
  {\small\textnormal{SysMA Research Unit, IMT School for Advanced
  Studies Lucca}}\\[0.8em]
  Joost J. Joosten\\[0.35em]
  {\small\textnormal{Department of Philosophy, University of
  Barcelona}}\\
  {\small\textnormal{and}}\\
  {\small\textnormal{Centre de Recerca Matemàtica, Barcelona}}}
\date{}

\makeatletter
\renewcommand{\maketitle}{%
  \begin{center}
    {\Large\scshape \@title\par}
    \ifx\@author\@empty\else
      \vspace{0.8em}{\itshape \@author\par}
    \fi
  \end{center}
  \vspace{1.1em}
}
\makeatother

\renewenvironment{abstract}
  {\begin{quote}\small\noindent\textsc{Abstract}\par\smallskip}
  {\end{quote}}

\newenvironment{keywords}
  {\par\smallskip\small\noindent\textit{Keywords.}\ }
  {\par}
\newenvironment{msc}
  {\small\noindent\textit{2020 Mathematics Subject Classification.}\ }
  {\par\medskip}

\begin{document}

\maketitle

\begin{abstract}
The reflection calculus with one modality is denoted by $\RC_1$ and its quantified generalization is known as $\QRC$. These decidable calculi are of interest since they generate the strictly positive fragments of the (quantified) provability logic for theories like Peano Arithmetic and Heyting Arithmetic (\cite{BorgesJoosten:2026:SPLHA,BorgesJoosten:2023:Escape}).

We introduce  analytic cut-free calculi \calc and \qcalc for  $\RC_1$ and $\QRC$ respectively. The new calculi derive sequents of the form $A\implies B$ with $A,B$ strictly positive formulas. In particular, the new calculi do not use (multi)sets nor nested sequents.

For both \calc and \qcalc it is shown that the cut-rule is admissible. As a consequence, it is proven that the calculi are deductively equivalent to $\RC_1$ and $\QRC$ respectively. The cutfree systems are used to prove uniform interpolation for the two reflection calculi. 

In addition, \calc is used to give a syntactic proof of Dashkov's Theorem to the extent that $\RC_1$ is the strictly positive fragment of the provability logic $\gl$. 

The paper also introduces non-analytic calculi that are conceptually easier but proof-theoretically behave badly. It is shown that even those systems allow direct cut-elimination: the benefits of the strict positivity outweigh the draw-backs of a non-analytic calculus in these cases.
\end{abstract}

\begin{msc}
\papermsc
\end{msc}

\section{Introduction}
\label{sec:introduction}

Strictly positive modal logics restrict the modal grammar to propositional variables, $\top$, conjunction, and diamond modalities. Their deductive presentations often concern expressions of shape $A\vdash B$ between formulas of this restricted language. In the field of provability logics~\cite{Boolos:1993:LogicOfProvability,sep-logic-provability}, strictly positive fragments retain a substantial value of the proof-theoretic content of the full modal systems and their arithmetical applications, in spite of their restricted language. In particular, reflection calculi provide a powerful framework for reasoning about iterated reflection principles and carrying out ordinal analysis \cite{Beklemishev:2012:Calibrating,Beklemishev:2014:PositiveReflection,
Beklemishev:2018:ConservativitySpectra,Beklemishev:2005:Survey,DBLP:conf/wollic/Beklemishev17}.

The unimodal reflection calculus $\RC_1$ corresponds to the strictly positive fragment of Gödel-Löb provability logic \gl~\cite{Dashkov:2012:PositiveFragment}. Its quantified extension \QRC, replaces propositional variables with relational atoms over a fixed signature, and extends the grammar by universal quantification, while retaining the $\RC_1$ operators~\cite{BorgesJoosten:2020:QRC1}. Besides its intrinsic value as a quantified modal logic, \QRC has a precise arithmetical interest: it captures the strictly positive fragment of the quantified provability logic for Peano Arithmetic and Heyting Arithmetic~\cite{BorgesJoosten:2023:Escape,BorgesJoosten:2026:SPLHA}.

The standard axiomatisations of reflection calculi, including $\RC_1$ and \QRC, include a transitivity rule for the relation $\vdash$, sometimes labelled ``syllogism'', which plays the role of Cut in standard sequent calculi. For proof analysis and automation of reasoning, though, it is generally preferable to work with systems where the cut rule can be effectively eliminated or proven to be admissible, instead of rooted in deductive engine of the logics.

Cut-free presentations of $\RC_1$ and \QRC could clearly be given in the context of enriched sequent calculi~\cite{DBLP:journals/rsl/Poggiolesi09,Negri:2005:ProofAnalysisModal,Shamkanov:2015:NestedGLP}. Recently, tree-rewriting calculi have been further developed for reflection calculi and related strictly positive logics, together with normalisation results~\cite{SantiagoFernandezEtAl:2024:TreeRC,
SantiagoFernandezEtAl:2025:TreeStrictPositive}. 

Here, we use none of those additional structures: we define $\RC_1$ and \QRC via Gentzen-style sequents consisting of one antecedent formula and one consequent formula, so that all rules of our calculi act on these formulas directly.

We first introduce an analytic calculus $\calc$ for $\RC_1$. Its modal rules combine necessitation/monotonicity \AxiomC{$A\vdash B$}\UnaryInfC{$\Diamond A\vdash\Diamond B$}\DisplayProof with diamond transitivity ($\Diamond\Diamond A\vdash\Diamond A$) in a form suited to backwards proof search: \AxiomC{$A\vdash \Diamond B$}\RightLabel{.}\UnaryInfC{$\Diamond A\vdash\Diamond B$}\DisplayProof The calculus $\qcalc$ for $\QRC$ is then obtained by passing to the quantified language and adjoining only the left and right rules for $\forall$. Next, we consider  alternative presentations (named \calcFourNec, \qcalcFourNec, and \calcFourNecStar, \qcalcFourNecStar) in which diamond transitivity is governed by explicit (left and) right rules, and show how they are related to the analytic calculi.

Our central technical result is the constructive proof of cut admissibility, in the shape of a cut elimination algorithm. For the propositional calculi, the reductions are organised by induction on the heights of the derivations above the proposed cut, without a separate induction on the complexity of the cut formula. The quantified setting follows the same pattern of cut-elimination.

Focussing on the analytic calculi, we show their relevance through two applications.
Firstly, we revisit Dashkov's characterisation of $\RC_1$ as the strictly positive fragment of $\gl$. Starting from a cut-free sequent proof for $\gl$, we extract a derivation in $\calc$, obtaining a novel and purely syntactic proof of the conservativity result.
Secondly, we study uniform interpolation in both reflection calculi. Craig interpolation for $\RC_1$ and \QRC was already observed in~\cite{BorgesJoosten:2026:SPLHA}; here we give a uniform
construction for the stronger task of forgetting a prescribed part of the vocabulary. In the quantified language, the construction removes relation symbols, individual constants and free individual parameters by replacing the relevant atomic occurrences with $\top$. Our constructions relate thus to the wider study of proof-theoretic techniques for interpolants~\cite{Iemhoff:2019:UniformInterpolation,DBLP:journals/corr/abs-2602-16318}. 

The results of this paper and its sequel have been sent as extended abstract on May 29, 2026 to \textit{The Proof Society Workshop on Proof Theory} where it has been accepted after peer-review and it was presented in Aussois on September 7, 2026. As it turns out, Nedeljkovic independently developed a cut-free calculus for $\QRC$ which is used to obtain similar results including NP-completeness and which he submitted June 2, 2026. His paper \cite{Nedeljkovic:2026:sequentcalculusquantifiedreflectionArxiv} is now online at the Arxiv too. His calculus works with finite multisets on the left instead of single formulas.

\medskip

The remaining part of the paper is organised as follows. Section~\ref{sec:preliminaries} fixes the languages and standard deductive systems for $\RC_1$ and \QRC. Section~\ref{sec:analytic-calculi} introduces our novel sequent calculi and establishes their basic properties; constructive cut admissibility and completeness are treated in Sections~\ref{sec:cut-n-compl} and~\ref{sec:cut-non-analytic}. Section~\ref{sec:gl-positive-fragment} gives the syntactic conservativity result over $\gl$-derivations, and Section~\ref{sec:modular-uniform-interpolation} proves the uniform interpolation theorem. 
We conclude, in~\Cref{sec:conclusion}, with questions and preliminary results about extensions to polymodal reflection calculi.

\section{Preliminaries}
\label{sec:preliminaries}

In this section we set the notation  used throughout the paper.

\noindent The strictly positive propositional language $\lang$ is generated by
$$A ::= \top \mid p \mid (A\land A) \mid \Diamond A$$ where $p$ comes from a countable set of propositional variables.
For the quantified language, we fix a countably infinite set
$\mathsf{Var}$ of individual variables, typically denoted by $x, y, \ldots$ and a first-order signature
$\Sigma$.  The signature contains relation symbols of fixed finite
arity and individual constants $c\in \mathsf{Cons}$, but no function symbols and no
equality.  Terms are therefore
$t ::= x \mid c$.
The strictly positive quantified language $\qlang$ is generated by
$$A ::= \top \mid S(t_1,\ldots,t_k) \mid (A\land A) \mid \forall x\,A \mid \Diamond A.$$
The quantifier-free fragment over nullary relation symbols can be identified
with $\lang$. 

The notation $t\subst{x}{t'}$ stands for substituting $t'$ for $x$ in $t$ and is usually defined by induction on $t$ which in our limited syntax reduces to the two base cases $t=c$ and $t=y$. Likewise, the notation $A\subst{x}{t}$ stands for substituting $t$ for $x$ in $A$ and is defined as usual by induction on $A$. We will almost exclusively be interested in cases where $t$ is free for $x$ in $A$ and sometimes explicitly require so. We refer to any textbook on basic proof theory like \cite{NegriVonPlato:2001:StructuralProofTheory} also for standard notions like the the multi-set of free variables and of bounded variables of a formula $A$, denoted $\FV(A)$ and $\BV(A)$ respectively.   With our limited syntax, the following observation becomes near immediate.
\begin{proposition}\label{theorem:NestedSubstitution}
    For $A$ a \qlang formula and for $t$ free for $x$ in $A$ and $t'$ free for $y$ in $A\subst{x}{t}$, with $x\neq y$ and $x\notin t'$, we have 
    \[
    \Big(A\subst{x}{t}\Big)\subst{y}{t'} = \Big(A\subst{y}{t'}\Big)\subst{x}{\big(t\subst{y}{t'}\big)}.
    \]
\end{proposition}
In this paper, we will give cut-free calculi for ${\RC_1}$ (\cite{Dashkov:2012:PositiveFragment}) and $\QRC$ (\cite{BorgesJoosten:2020:QRC1}) whose definitions we recall here.

\begin{definition}
The consequence relation $\vdash_{\RC_1}$ is the least relation on $\lang$
closed under the schemes in \Cref{fig:rc1-deductive-rules}.  

\begin{figure}[H]
\centering
\begin{tabular}{c@{\qquad}c@{\qquad}c}
  \AxiomC{}
  \RightLabel{\scriptsize id}
  \UnaryInfC{$A\vdash_{\RC_1}A$}
  \DisplayProof
&
  \AxiomC{}
  \RightLabel{\scriptsize $\top$}
  \UnaryInfC{$A\vdash_{\RC_1}\top$}
  \DisplayProof
&
  \AxiomC{$A\vdash_{\RC_1}B$}
  \AxiomC{$B\vdash_{\RC_1}C$}
  \RightLabel{\scriptsize Cut}
  \BinaryInfC{$A\vdash_{\RC_1}C$}
  \DisplayProof
\\[3ex]
  \AxiomC{}
  \RightLabel{\scriptsize $\land\mathcal L_1$}
  \UnaryInfC{$A\land B\vdash_{\RC_1}A$}
  \DisplayProof
&
  \AxiomC{}
  \RightLabel{\scriptsize $\land\mathcal L_2$}
  \UnaryInfC{$A\land B\vdash_{\RC_1}B$}
  \DisplayProof
&
  \AxiomC{$A\vdash_{\RC_1}B$}
  \AxiomC{$A\vdash_{\RC_1}C$}
  \RightLabel{\scriptsize $\land\mathcal R$}
  \BinaryInfC{$A\vdash_{\RC_1}B\land C$}
  \DisplayProof
\\[3ex]
  \AxiomC{$A\vdash_{\RC_1}B$}
  \RightLabel{\scriptsize Nec}
  \UnaryInfC{$\Diamond A\vdash_{\RC_1}\Diamond B$}
  \DisplayProof
&
  \multicolumn{2}{c}{%
    \AxiomC{}
    \RightLabel{\scriptsize 4}
    \UnaryInfC{$\Diamond\Diamond A\vdash_{\RC_1}\Diamond A$}
    \DisplayProof}
\end{tabular}
\caption{Deductive presentation of the unimodal reflection calculus
$\RC_1$.}
\label{fig:rc1-deductive-rules}
\end{figure}
The quantified reflection calculus $\QRC$ is the least extension of the
same positive core to $\qlang$ that is also closed under the four schemes in
\Cref{fig:qrc-deductive-rules}.  
\begin{figure}[H]
\centering
\begin{tabular}{c@{\qquad}c}
  \AxiomC{$A\vdash_{\QRC}B$}
  \RightLabel{\scriptsize $\forall\mathcal R$}
  \UnaryInfC{$A\vdash_{\QRC}\forall x\,B$}
  \DisplayProof
&
  \AxiomC{$A\subst{x}{t}\vdash_{\QRC}B$}
  \RightLabel{\scriptsize $\forall\mathcal L$}
  \UnaryInfC{$\forall x\,A\vdash_{\QRC}B$}
  \DisplayProof
\\[1.2ex]
  {\scriptsize $x\notin\FV(A)$}
&
  {\scriptsize $t$ free for $x$ in $A$}
\\[3.2ex]
  \AxiomC{$A\vdash_{\QRC}B$}
  \RightLabel{\scriptsize Inst}
  \UnaryInfC{$A\subst{x}{t}\vdash_{\QRC}B\subst{x}{t}$}
  \DisplayProof
&
  \AxiomC{$A\subst{x}{c}\vdash_{\QRC}B\subst{x}{c}$}
  \RightLabel{\scriptsize CE}
  \UnaryInfC{$A\vdash_{\QRC}B$}
  \DisplayProof
\\[1.2ex]
  {\scriptsize $t$ free for $x$ in both $A$ and $B$}
&
  {\scriptsize $c$ occurs in neither $A$ nor $B$}
\end{tabular}
\caption{The additional rules of the deductive calculus $\QRC$.}
\label{fig:qrc-deductive-rules}
\end{figure}

\end{definition}

We observe that no additional interaction rule between $\Diamond$ and $\forall$ is primitive
 in this standard presentation.\footnote{In particular, the useful principle
$\Diamond\forall x\,A\vdash_{\QRC}\forall x\,\Diamond A$ is derivable from
the displayed rules.} We fix some standard metrics on our formulas: modal depth $\mdepth$ and quantifier depth $\qdepth$. Likewise, formula size $|A|$ counts logical symbols, atomic relation symbols, and term occurrences:
\ \\
\begin{minipage}[t]{0.6\textwidth}
\[
\begin{array}{c|cc}
  C & \mdepth(C) & \qdepth(C) \\
  \hline
  \top, p, \ S(\vec t) & 0 & 0 \\
  A\land B
    & \max\{\mdepth(A),\mdepth(B)\}
    & \max\{\qdepth(A),\qdepth(B)\} \\
  \forall x\,A & \mdepth(A) & \qdepth(A)+1 \\
  \Diamond A & \mdepth(A)+1 & \qdepth(A).
\end{array}
\]
\end{minipage}
\begin{minipage}[t]{0.4\textwidth}
\[
\begin{aligned}
  |\top| &= 1,\\
  |S(t_1,\ldots,t_k)| &= k+1,\\
  |A\land B| &= |A|+|B|+1,\\
  |\forall x\,A| &= |A|+1,\\
  |\Diamond A| &= |A|+1.
\end{aligned}
\]
\end{minipage}

We see that terms are either variables or constants, and  substitution of a term for a variable  preserves logical shape,
formula size, modal depth, and quantifier depth.

\section{Sequent calculi design and analysis}
\label{sec:analytic-calculi}

All Gentzen systems in this section use sequents of the form
$A\implies B$ where $A$ and $B$ are single formulas. In particular, we will not use
sets or multisets of formulas as contexts \cite{TroelstraSchwichtenberg:2000:BasicProofTheory}, we do not use labels \cite{Negri:2005:ProofAnalysisModal}, or nesting of multisets inside multisets as in \cite{Shamkanov:2015:NestedGLP}.  

If $\Pi$ is a derivation, a
zero-premise inference has height $0$, while an inference with immediate
subderivations $\Pi_1,\ldots,\Pi_k$ has height
$1+\max_i h(\Pi_i)$. For a calculus $\mathsf G$ we write
$\mathsf G\vdash_{\leq n}A\implies B$
when there exists a derivation $\Pi$ of height at most $n$.  Exact height is
written only for a fixed derivation as $h(\Pi)=n$. Sometimes we will write $\mathcal D \vdash_{\mathsf G} A\implies B$ to fix a particular $\mathsf G$ derivation $\mathcal{D}$ of the sequent $A\implies B$ which is called the conclusion of $\mathcal{D}$. 
We will sometimes write 
\AxiomC{$\mathcal{D}$}\noLine\UnaryInfC{$A\implies B$}\DisplayProof to denote a proof whose conclusion (last/bottom line) is $A \implies B$ so that in particular $\mathcal D$'s height can be null under this writing convention.

Besides maintaining sequents as pairs of formulas, we design our sequent calculi as modularly as the original reflection calculi: the lifting to first-order syntax of \qlang involves just adding a concise module of rules to the propositional core, for handling the universal quantifier naturally.

\subsection{Analytic systems and soundness}
 \label{subsec:analytic-core}

\begin{definition}[The calculus \calc]
\label{def:modular-g1rc1}
The sequents of \calc are $A\implies B$, with $A,B\in\lang$. Its deductive rules are those in \Cref{fig:g1rc1-rules}.
\end{definition}

\begin{figure}[H]
\centering
\begin{tabular}{c@{\qquad\qquad}c}
  \AxiomC{}
  \RightLabel{\scriptsize Initial}
  \UnaryInfC{$p\implies p$}
  \DisplayProof
&
  \AxiomC{}
  \RightLabel{\scriptsize $\top\mathcal R$}
  \UnaryInfC{$A\implies\top$}
  \DisplayProof
\\[3ex]
  \AxiomC{$A\implies C$}
  \RightLabel{\scriptsize $\land\mathcal L_1$}
  \UnaryInfC{$A\land B\implies C$}
  \DisplayProof
&
  \AxiomC{$B\implies C$}
  \RightLabel{\scriptsize $\land\mathcal L_2$}
  \UnaryInfC{$A\land B\implies C$}
  \DisplayProof
\\[3ex]
  \AxiomC{$A\implies B$}
  \AxiomC{$A\implies C$}
  \RightLabel{\scriptsize $\land\mathcal R$}
  \BinaryInfC{$A\implies B\land C$}
  \DisplayProof
&
  \AxiomC{$A\implies B$}
  \RightLabel{\scriptsize Nec}
  \UnaryInfC{$\Diamond A\implies\Diamond B$}
  \DisplayProof
\\[3ex]
  \multicolumn{2}{c}{%
    \AxiomC{$A\implies\Diamond B$}
    \RightLabel{\scriptsize $4\mathrm{Nec}$}
    \UnaryInfC{$\Diamond A\implies\Diamond B$}
    \DisplayProof}
\end{tabular}
\caption{Rules of the analytic calculus \calc.}
\label{fig:g1rc1-rules}
\end{figure}
The quantified analytic calculus \qcalc is obtained by adding a left and a right rule for universal quantification to \calc.

\begin{definition}[The calculus \qcalc]
\label{def:modular-g1qrc1}
The calculus \qcalc is obtained from \calc by changing the
language from \lang to \qlang
and adding the two quantifier rules {\scriptsize $\forall\mathcal L$} and {\scriptsize $\forall\mathcal R$} as displayed in
\Cref{fig:g1qrc1-quantifier-rules}.  All other rules are unchanged, except for initial sequents having shape $R(\vec{t})\implies R(\vec{t})$ for any relation $R$ in the language.
\end{definition}

\begin{figure}[H]
\centering
\begin{tabular}{c@{\qquad\qquad}c}
  \AxiomC{$A\subst{x}{t}\implies B$}
  \RightLabel{\scriptsize $\forall\mathcal L$}
  \UnaryInfC{$\forall x\,A\implies B$}
  \DisplayProof
&
  \AxiomC{$A\implies B\subst{x}{y}$}
  \RightLabel{\scriptsize $\forall\mathcal R$}
  \UnaryInfC{$A\implies\forall x\,B$}
  \DisplayProof
\\[1.2ex]
  {\scriptsize $t$ free for $x$ in $A$}
&
  {\scriptsize $y$ free for $x$ in $B$; $y\notin\FV(A)\cup\FV(\forall x\,B)$}
\end{tabular}
\caption{Quantifier rules added in the analytic calculus \qcalc.}
\label{fig:g1qrc1-quantifier-rules}
\end{figure}

The variable $y$ in $\forall\mathcal R$ is the eigenvariable of the inference; it is displayed as a separate metavariable from the bound variable $x$.  Thus, the side condition for $\forall\mathcal R$ says that $y$ is not free in the conclusion in addition to demanding $y$ being free for $x$ in $B$.  
This is the standard G3-style formulation used in the proof-theoretic treatment of bound-variable conversion and substitution~\cite[Section~4.1]{NegriVonPlato:2001:StructuralProofTheory}.

\begin{lemma}[Generalised initial sequents]
\label{lem:modular-generalised-id}
For every $A\in\lang$,
$\calc\vdash A\implies A$.
Likewise, for every $A\in\qlang$,
$\qcalc\vdash A\implies A$.
\end{lemma}

\begin{proof}
Both claims are proved by induction on $A$.  The cases for $\top$,
atoms, conjunction, and $\Diamond$ are simple and identical in the two calculi.
Thus we focus on quantified formulas with $A=\forall x\,B$.  
Choose $y$ fresh for $B$; by the induction hypothesis,
$B\subst{x}{y}\implies B\subst{x}{y}$ is derivable.  Hence,
\begin{prooftree}
  \AxiomC{$B\subst{x}{y}\implies B\subst{x}{y}$}
  \RightLabel{\scriptsize $\forall\mathcal L$}
  \UnaryInfC{$\forall x\,B\implies B\subst{x}{y}$}
  \RightLabel{\scriptsize $\forall\mathcal R$}
  \UnaryInfC{$\forall x\,B\implies\forall x\,B$}
\end{prooftree}
where the same fresh $y$ is used as the instance term and as the eigenvariable of the final rule.
\end{proof}
It is easy to see that the two analytic calculi \calc and \qcalc do not prove more theorems than the original calculi $\RC_1$ and $\QRC$ respectively.
\begin{theorem}[Soundness]
\label{thm:modular-soundness}
For $A,B\in\lang$, if $\calc\vdash A\implies B$, then
$A\vdash_{\RC_1}B$.  Likewise, for $A,B\in\qlang$, if
$\qcalc\vdash A\implies B$, then $A\vdash_{\QRC}B$.
\end{theorem}

\begin{proof}
We induct on the given derivation in \calc or \qcalc.  The common rules are interpreted
in the same way in the two deductive systems. In particular, the $\wedge \mathcal L_i$ rule in the sequent calculi is taken care of by the $\wedge \mathcal L_i$ axiom, and a cut on the inductive hypothesis. For $4\mathrm{Nec}$, from the inductive hypothesis $A\vdash\Diamond B$ we obtain
$\Diamond A\vdash\Diamond\Diamond B$ by necessitation and then use
modal transitivity $\Diamond\Diamond B\vdash\Diamond B$ and Cut.
In the quantified calculus, $\forall\mathcal L$ is exactly the
corresponding clause of the deductive presentation of \QRC. For $\forall\mathcal R$, suppose the last inference is
\begin{prooftree}
  \AxiomC{$A\implies B\subst{x}{y}$}
  \RightLabel{\scriptsize $\forall\mathcal R$ \ .}
  \UnaryInfC{$A\implies\forall x\,B$}
\end{prooftree}
By the induction hypothesis,
$A\vdash_{\QRC}B\subst{x}{y}$.  Let $c, d$ be a fresh constants. We make a distinction between the cases $x\notin \FV (A)$ and $x\in \FV (A)$. In case $x\notin \FV (A)$ we obtain a variable renaming in $\QRC$ as follows:
\begin{prooftree}
  \AxiomC{$A\implies B\subst{x}{y}$}
  \RightLabel{\scriptsize Inst}
  \UnaryInfC{$A\implies (B\subst{x}{y})\subst{y}{c}$}
  \RightLabel{\scriptsize CE}
  \UnaryInfC{$A\implies \Big((B\subst{x}{y})\subst{y}{c}\Big) \subst{c}{x}$}
  \RightLabel{\scriptsize $\forall\mathcal R$ \ .}
  \UnaryInfC{$A\implies\forall x\,B$}
\end{prooftree}
We observe that all the side conditions are met since $x\notin \FV (A)$ and $y \notin \FV (A)$.

In case $x\in \FV (A)$ we obtain two variable renamings in $\QRC$ as follows:
\begin{prooftree}
  \AxiomC{$A\implies B\subst{x}{y}$}
  \RightLabel{\scriptsize Inst}
  \UnaryInfC{$A\implies (B\subst{x}{y})\subst{y}{c}$}
  \RightLabel{\scriptsize Inst}
  \UnaryInfC{$A\subst{x}{d}\implies (B\subst{x}{y})\subst{y}{c} $}

  \RightLabel{\scriptsize CE}
  \UnaryInfC{$A\subst{x}{d}\implies \Big((B\subst{x}{y})\subst{y}{c}\Big) \subst{c}{x}$}
  \RightLabel{\scriptsize $\forall\mathcal R$}
  \UnaryInfC{$A\subst{x}{d}\implies\forall x\,B$}
  \RightLabel{\scriptsize CE \ \ 
  \ \ \ \ \ \ \ .}
  \UnaryInfC{$A\implies\forall x\,B$}
\end{prooftree}
Again, it is easy to see that all the side-conditions are satisfied.

\end{proof}

\subsection{Non-analytic systems and equivalence}
\label{sec:modular-cut}

In this section, we present two additional non-analytic sequent systems for the propositional reflection calculus, together with their quantified extensions.  They are of specific structural interest: the first system (\calcFourNec) exposes diamond transitivity by separate left and right rules, whereas the second one (\calcFourNecStar) shows that the left transitivity rule is ultimately redundant.  Albeit their non-analiticity, we shall see that both admit a satisfactory proof analysis.

\begin{definition}[The calculi \calcFourNec and \qcalcFourNec]
\label{def:g1k4-plus}
The calculus \calcFourNec has Initial, $\top\mathcal R$,
$\land\mathcal L_i$, $\land\mathcal R$, and Nec from \calc (Figure~\ref{fig:g1rc1-rules}), while
$4\mathrm{Nec}$ is replaced by the two rules displayed in
\Cref{fig:g1k4-plus-modal-rules}. The calculus \qcalcFourNec is obtained from \calcFourNec by enlarging
the language to \qlang, allowing arbitrary relational atoms in Initial,
and adjoining $\forall\mathcal L$ and $\forall\mathcal R$ exactly as in
\qcalc (Figure~\ref{fig:g1qrc1-quantifier-rules}).
\end{definition}

\begin{figure}[H]
\centering
\begin{tabular}{c@{\qquad\qquad}c}
  \AxiomC{$\Diamond A\implies B$}
  \RightLabel{\scriptsize $4\mathcal L$}
  \UnaryInfC{$\Diamond\Diamond A\implies B$}
  \DisplayProof
&
  \AxiomC{$A\implies\Diamond\Diamond B$}
  \RightLabel{\scriptsize $4\mathcal R$}
  \UnaryInfC{$A\implies\Diamond B$}
  \DisplayProof
\end{tabular}
\caption{The two explicit transitivity rules of \calcFourNec.}
\label{fig:g1k4-plus-modal-rules}
\end{figure}
Thus, the passage from \calcFourNec to \qcalcFourNec is again a pure
quantifier module. We observe that the $4 \mathcal{L}$ rule is analytic whereas the $4 \mathcal{R}$ is not. We now introduce the starred systems where we only keep the non-analytic extra rule.

\begin{definition}[The calculi \calcFourNecStar and \qcalcFourNecStar]
\label{def:g1k4-plus-star}
The calculus \calcFourNecStar is obtained from \calcFourNec by deleting
$4\mathcal L$.  Likewise, \qcalcFourNecStar is obtained from
\qcalcFourNec by deleting $4\mathcal L$.
\end{definition}

The star calculi therefore retain only the non-analytic right rule
$4\mathcal R$ for diamond transitivity. The following observation is immediate.

\begin{proposition}\label{Theorem:StarPlusContainedInPlus}
    $\calcFourNecStar \subseteq \calcFourNec$ and, $\qcalcFourNecStar \subseteq \qcalcFourNec$.
\end{proposition}
The only rule from \calc that is not present in $\calcFourNecStar$ is the $4\mathrm{Nec}$ rule. It is easy to see that this rule is derivable in \calcFourNecStar.
\begin{lemma}[Derivability of $4\mathrm{Nec}$]
\label{lem:four-nec-in-reduced}
The rule
\begin{prooftree}
  \AxiomC{$A\implies\Diamond B$}
  \RightLabel{\scriptsize $4\mathrm{Nec}$}
  \UnaryInfC{$\Diamond A\implies\Diamond B$}
\end{prooftree}
is derivable in \calcFourNec and \calcFourNecStar, and likewise in their
quantified counterparts.
\end{lemma}

\begin{proof}
Use Nec followed by $4\mathcal R$:
\begin{prooftree}
  \AxiomC{$A\implies\Diamond B$}
  \RightLabel{\scriptsize Nec}
  \UnaryInfC{$\Diamond A\implies\Diamond\Diamond B$}
  \RightLabel{\scriptsize $4\mathcal R$ \ \ .}
  \UnaryInfC{$\Diamond A\implies\Diamond B$}
\end{prooftree}
\end{proof}

As a corollary, we get the following extension of Proposition \ref{Theorem:StarPlusContainedInPlus}.

\begin{proposition}\label{Theorem:CalcContainedInStarPlusContainedInPlus}
    $\calc \subseteq \calcFourNecStar \subseteq \calcFourNec$ and, $\qcalc \subseteq \qcalcFourNecStar \subseteq \qcalcFourNec$.
\end{proposition}

\begin{proof}
    Building upon Proposition \ref{Theorem:CalcContainedInStarPlusContainedInPlus} we see by induction on the length of a \calc derivation that $\calc \subseteq \calcFourNecStar$. The only non-trivial case is the 4Nec rule which is covered by Lemma \ref{lem:four-nec-in-reduced}. The quantifier version directly follows suit.
\end{proof}

We shall now see that in this Proposition \ref{Theorem:CalcContainedInStarPlusContainedInPlus}, we actually have equalities and not just inclusions. To achieve this, it suffices to prove admissibility of both the $4\mathcal R$ and the $4\mathcal L$ rules in \calc.

\begin{lemma}[Height-preserving admissibility of $4\mathcal R$]
\label{lem:modular-4r-admissible}
Let $\mathsf G$ be either \calc or \qcalc.\\
If
$\mathsf G\vdash_{\leq n}A\implies\Diamond\Diamond B$, then
$\mathsf G\vdash_{\leq n}A\implies\Diamond B$.
\end{lemma}

\begin{proof}
By an easy induction on the derivation of $A\implies\Diamond\Diamond B$.  We first notice that zero-premise rules and $\land\mathcal R$ cannot be the last applied rule;
in the quantified calculus
neither can $\forall\mathcal R$.  

If the last inference is a left rule, we apply the induction hypothesis (only one indeed) to its premise and then apply the left rule again to get the result.  If the last inference is Nec, its premise is of the form $A_0\implies\Diamond B$, so that one application of
$4\mathrm{Nec}$ gives the required conclusion.  If the last inference is $4\mathrm{Nec}$, we first apply the induction hypothesis to its premise and
then apply $4\mathrm{Nec}$ again.  

It is easy to see that heights are always preserved.  
\end{proof}
We have height preserving admissibility of $4\mathcal R$. The admissibility of $4\mathcal L$ is proven without preservation of height since we have to include a single additional step in the proof.
\begin{lemma}[Admissibility of $4\mathcal L$]
\label{lem:modular-4l-admissible}
Let $\mathsf G$ be either \calc or \qcalc. \\
If
$\mathsf G\vdash_{\leq n}\Diamond A\implies C$, then
$\mathsf G\vdash_{\leq n{+}1}\Diamond\Diamond A\implies C$.
\end{lemma}

\begin{proof}

Assume that $\mathsf G\vdash\Diamond A\implies C$. Given the form of this sequent it can only admit a zero-height proof or one where the last rule was either $\wedge \mathcal{R}$, $\mathrm{Nec}$, $4\mathrm{Nec}$ or, in the case of \qcalc, $\forall \mathcal{R}$. The zero-height case can only be with $C=\top$, where the claim is immediate.  

If the final rule is either $\mathrm{Nec}$ or $4\mathrm{Nec}$, the conclusion is of the form $\Diamond A\implies\Diamond C_0$ so that we apply $4\mathrm{Nec}$ to conclude $\Diamond\Diamond A\implies\Diamond C_0$. We observe that these cases do not invoke the inductive hypothesis and one adds a single step to the proof.

If the last rule was either $\wedge \mathcal{R}$ or $\forall \mathcal R$ we apply the inductive hypothesis to the antencedent(s) and apply the rule again. For example, in the quantified calculus, if $\forall \mathcal R$ was the last rule we have that $C=\forall x\,C_0$. We apply the inductive hypothesis to $\Diamond A\implies C_0\subst{x}{y}$ to obtain $\Diamond\Diamond A\implies C_0\subst{x}{y}$ and we get the result by applying $\forall \mathcal{R}$ again (observing that the side conditions are maintained). The $\leq n{+}1$ condition is easy to check.

\end{proof}

We can now prove the equivalence of the three calculi $\calc, \calcFourNec$ and $\calcFourNecStar$ and likewise for their quantified versions.

\begin{theorem}[Equivalence of sequent presentations]
\label{thm:modular-gentzen-equivalence}
For $A,B\in\lang$, the following equivalences hold:
\[
\calc\vdash A\implies B
\quad\textit{iff}\quad
\calcFourNec\vdash A\implies B
\quad\textit{iff}\quad
\calcFourNecStar\vdash A\implies B.
\]

For $A,B\in\qlang$, the following equivalences hold:
\[
\qcalc\vdash A\implies B
\quad\textit{iff}\quad
\qcalcFourNec\vdash A\implies B
\quad\textit{iff}\quad
\qcalcFourNecStar\vdash A\implies B.
\]

Moreover, all implications are effective transformations of derivations in the calculi under analysis.
\end{theorem}

\begin{proof}
In Propositon \ref{Theorem:CalcContainedInStarPlusContainedInPlus} we already proved that     $\calc \subseteq \calcFourNecStar \subseteq \calcFourNec$ and, $\qcalc \subseteq \qcalcFourNecStar \subseteq \qcalcFourNec$ so that we only need to show that $\calcFourNec \subseteq \calc$ and that $\qcalcFourNec \subseteq \qcalc$. This follows via an easy induction on the length of a proof where the only interesting steps are $4 \mathcal R$ and $4 \mathcal L$ which are taken care of by Lemma \ref{lem:modular-4r-admissible} and Lemma \ref{lem:modular-4l-admissible} respectively.
\end{proof}

\section{Cut elimination and completeness for the analytic calculi}
\label{sec:cut-n-compl}
The Cut rule
\begin{prooftree}
    \AxiomC{$A\implies C$}\AxiomC{$C\implies B$}
\RightLabel{\textit{Cut}}\BinaryInfC{$A\implies B$}
\end{prooftree}
is not primitive in any of the calculi \calc, \calcFourNec, \calcFourNecStar nor in any of their quantified extensions.  In this section we will prove the admissibility of the cut-rule for any of our systems, \textit{constructively}. The easiest proof is for \calc and \qcalc and admissibility for the other systems easily follows. 

It is well known that strictly positive logics are often tame and behave better than their full counterparts. In the case of proof-calculi we will showcase this phenomenon once more.\footnote{Besides for \calc, in~\Cref{sec:cut-non-analytic}, we will give a direct proof of cut elimination for \calcFourNec: even though this is a non-analytic and unelegant calculus (technically speaking), the tameness of the strictly positive fragment allows a fairly simple direct proof of cut-elimination.}

\subsection{Cut elimination for \calc}

For our proof analysis, we refer to the standard proof theoretic terminology: we say that
a formula is \emph{principal} in an inference when its outermost connective
is introduced by that inference.  For the \calc calculus we say that $A\land B$ is principal in
$\land\mathcal L_i$, the consequent $B\land C$ in
$\land\mathcal R$, both $\Diamond A$ and $\Diamond B$ are principal in {\scriptsize{Nec}} and {\scriptsize{4Nec}}.   In the quantified systems, $\forall x\,A$ is
principal in $\forall\mathcal L$ when it occurs in the antecedent and in
$\forall\mathcal R$ when it occurs in the succedent.  The two occurrences
of the atom in {Initial} are principal, as is $\top$ in $\top\mathcal R$.\footnote{Refer to, e.g.,~\cite{TroelstraSchwichtenberg:2000:BasicProofTheory} and~\cite{NegriVonPlato:2001:StructuralProofTheory}.}

\begin{theorem}[Cut admissibility for \calc]
\label{thm:calc-cut}
If
$\calc\vdash A\implies B$ and
$\calc\vdash B\implies C$, then
$\calc\vdash A\implies C$.
\end{theorem}
\begin{proof}
    We perform an induction on the \textit{level} of the cut which is defined as the sum of the heights of the derivations for $A\implies B$ and $B\implies C$. 

    We make a case distinction: either $B$ is not principal in at least one of the derivations, or $B$ is principal in both derivations. We first see that the level can only be null when $B$ is principal in both proofs. This can only happen when $B=C=\top$ or when $A=B=C=p$ for some propositional variable $p$. In both cases we see that $A\implies C$ is again a derivable with height zero. Thus, our case distinction includes that base case of the induction.

We first consider the case when $B$ is not active in one of the derivations. We focus on the derivation $\mathcal D$ of the sequent where $B$ is not principal. We call the other proof $\mathcal D'$. If the height of that derivation is zero, this can only happen for $B\implies C$ with $C=\top$ in which case $A\implies C$ is evidently provable. So, the remaining situation is when $\mathcal D$ follows from premisses (possibly just one) with proof of smaller height. Each premiss contains $B$ so, by the induction hypothesis, we can cut the $B$ away since the level of the cuts has decreased.

We now consider the case that $B$ is active in both proofs. In the above, we already considered the easy cases where $B=\top$ or $B=p$. There are two remaining subcases.

When $B= B_0\wedge B_1$, we have derivations of {$A\implies B_0$}, {$A\implies B_1$} and $B_i\implies C$ for some $i\in \{0,1\}$. We can apply the inductive hypothesis since the cut on $A\implies B_i$ and $B_i \implies C$ to obtain $A\implies C$ is clearly of lower level.

When $B= \Diamond B_0$, on both premisses of the Cut, the last rule can have been either an application of {Nec} or of {4Nec} leaving us with the four scenarios from~\Cref{tab:placeholder}.

\begin{table}[h]
    \centering
    \begin{tabular}{|c|c|c|}
    \hline
        Scenario number & Last step in $A\implies \Diamond B_0$ & Last step in $\Diamond B_0 \implies C$\\
        \hline
         {\MakeUppercase{\romannumeral 1}}& {\scriptsize 4Nec}  & {\scriptsize 4Nec}\\
         \hline
         {\MakeUppercase{\romannumeral 2}}& {\scriptsize 4Nec}  & {\scriptsize Nec}\\
         \hline
         {\MakeUppercase{\romannumeral 3}}& {\scriptsize Nec}  & {\scriptsize 4Nec}\\
         \hline
         {\MakeUppercase{\romannumeral 4}}& {\scriptsize Nec}  & {\scriptsize Nec}\\
         \hline
    \end{tabular}
    \caption{Four scenarios for a cut with $\Diamond B_0$ principal for both premisses.}
    \label{tab:placeholder}
\end{table}
    In all scenarios we have that $A=\Diamond A_0$ and $C=\Diamond C_0$ for some $A_0, C_0$. The first two scenarios are treated the same: 
    \begin{prooftree}
    \AxiomC{$A_0\implies \Diamond B_0$}\AxiomC{$\Diamond B_0\implies \Diamond C_0$}
\RightLabel{\textit{Cut}  \ \ \  .}\BinaryInfC{$A_0\implies \Diamond C_0$}
\RightLabel{{\scriptsize 4Nec}}\UnaryInfC{$\Diamond A_0\implies \Diamond C_0$}
\end{prooftree}
The inductive hypothesis can be applied since the level of the cut on $\Diamond B_0$ is of lower level as the height of the proof of $A_0 \implies \Diamond B_0$ decreased. Scenarios {\MakeUppercase{\romannumeral 3}} and {\MakeUppercase{\romannumeral 4}} follow the same structure: cut with lower level the $B_0$ of both antecedents and then apply the same right rule again:  

\begin{minipage}[t]{0.5\textwidth}
\begin{prooftree}
    \AxiomC{$A_0\implies  B_0$}\AxiomC{$ B_0\implies \Diamond C_0$}
\RightLabel{\textit{Cut} }\BinaryInfC{$A_0\implies \Diamond C_0$}
\RightLabel{{\scriptsize 4Nec}}\UnaryInfC{$\Diamond A_0\implies \Diamond C_0$}
\end{prooftree}
\begin{center}
Scenario {\MakeUppercase{\romannumeral 3}}
\end{center}
\end{minipage}
\begin{minipage}[t]{0.5\textwidth}
\begin{prooftree}
    \AxiomC{$A_0\implies  B_0$}\AxiomC{$ B_0\implies  C_0$}
\RightLabel{\textit{Cut}  }\BinaryInfC{$A_0\implies  C_0$}
\RightLabel{{\scriptsize Nec}}\UnaryInfC{$\Diamond A_0\implies \Diamond C_0$}
\end{prooftree}
\begin{center}
Scenario {\MakeUppercase{\romannumeral 4}} \ \ \ \ .
\end{center}
\end{minipage}

\end{proof}

We stress that there is no induction on the cut-rank nor any nested induction in our proof. 

\subsection{Cut elimination with quantifiers}
We will now extend our result to \qcalc. It is easy to see that $\qcalc \nvdash \forall x\, \forall y \ \big(P(x) \wedge Q(y) \big) \implies P(y)$ yet it can be proven using the \textit{Cut} rule since both $\qcalc \vdash \forall x\, \forall y \ \big(P(x) \wedge Q(y) \big) \implies \forall z\  P(z)$ and $\qcalc \vdash \forall z\  P(z) \implies P(y)$.\footnote{The example is directly taken from \cite{TroelstraSchwichtenberg:2000:BasicProofTheory}.} Thus, \emph{strictly speaking}, \qcalc does not admit Cut-elimination.

Clearly, in the example this non-admissibility of Cut can be avoided if we had chosen the bounded variables differently.
Probably, the conceptually easiest way to deal with these binding issues is by demanding bound and free variables sui generis to come from disjoint sets. Curiously, this is not what is mostly done in the literature and one typically only requires that locally bound and free variables do not coincide and rename bound variables as needed. 

\begin{definition}
    For a \qcalc derivation $\mathcal{D}$ we define the multi-set of \emph{eigenvariables} ${\mathsf{ev}}(\mathcal{D})$ and the multi-set of \emph{bound variables} ${\mathsf{bv}}(\mathcal{D})$ inductively. The eigenvariables are the variables in a proof that disappear at some point due to a universal quantifier introduction so that ${\mathsf{ev}}(\mathcal{D}) :=  \varnothing$ when $\mathcal{D}$ has height zero. The bounded variables of $\mathcal{D}$ are the variables that occur bounded somewhere in $\mathcal{D}$, so that $ {\mathsf{bv}}(A\implies B) := \BV(A)$ when the derivation $\mathcal{D}$ of $A\implies B$ has height zero.
    \[
    \mathsf{ev}
    \left(\AxiomC{$\mathcal D$}
    \noLine
    \UnaryInfC{$A\implies B\subst{x}{y}$}
    \RightLabel{\scriptsize$\forall\mathcal  R$}\UnaryInfC{$A\implies \forall x B$} \DisplayProof
    \right) := {\mathsf{ev}}\left(\AxiomC{$\mathcal D$}
    \noLine
    \UnaryInfC{$A\implies B\subst{x}{y}$}
     \DisplayProof
    \right) \cup \{ y\}
    \]
    where $\cup$ denotes multiset union, and
    \[
    \mathsf{ev}\left(
      \AxiomC{$\mathcal D$}
      \RightLabel{\scriptsize $\mathcal R$}
      \UnaryInfC{$A \implies B$}
      \DisplayProof
    \right)
    := \mathsf{ev}(\mathcal D)
    \]
    for any other unary rule $\mathcal R$. 
    Similarly,  
    \[
    {\mathsf{bv}}\left(\AxiomC{$\mathcal D$}
    \noLine
    \UnaryInfC{$A\implies B\subst{x}{y}$}
    \RightLabel{\scriptsize$\forall\mathcal  R$}\UnaryInfC{$A\implies \forall x B$} \DisplayProof
    \right) := {\mathsf{bv}}\left(\AxiomC{$\mathcal D$}
    \noLine
    \UnaryInfC{$A\implies B\subst{x}{y}$}
     \DisplayProof
    \right) \cup \{ x\}
    \]
and
\[
    {\mathsf{bv}}\left(\AxiomC{$\mathcal D$}
    \noLine
    \UnaryInfC{$A\subst{x}{t}\implies B$}
    \RightLabel{\scriptsize$\forall\mathcal  L$}\UnaryInfC{$\forall x A\implies  B$} \DisplayProof
    \right) := {\mathsf{bv}}\left(\AxiomC{$\mathcal D$}
    \noLine
    \UnaryInfC{$A\subst{x}{t}\implies B$}
     \DisplayProof
    \right) \cup \{ x\} \ .
    \]
    
    For conjunction we state  
    \[
    {\mathsf{bv}}\left(\AxiomC{$\mathcal D$}\noLine\UnaryInfC{$A \implies C$}\RightLabel{\scriptsize$\wedge\mathcal{L}_1$} \UnaryInfC{$A \wedge B \implies C$} \DisplayProof\right) := {\mathsf{bv}}\left(\AxiomC{$\mathcal D$}\noLine\UnaryInfC{$A \implies C$}\RightLabel{\scriptsize$\wedge\mathcal{L}_2$} \UnaryInfC{$B \wedge A \implies C$} \DisplayProof\right) := {\mathsf{bv}}\left(\AxiomC{$\mathcal{D}$}\noLine \UnaryInfC{$A\implies C$}\DisplayProof\right) \cup \BV(B)
    \]
    and ${\mathsf{bv}}\left(\AxiomC{$\mathcal{D}$}\noLine\UnaryInfC{$A\implies B$}\DisplayProof\right) := {\mathsf{bv}} (\mathcal{D})$ for any other unary rule $\mathcal R$. Finally, the only binary rule is $\wedge \mathcal R$ and, the eigenvariables of a proof ending on $\wedge \mathcal R$ is defined to be the multi-set union of the eigenvariables of the two derivations leading to the antecedents of the rule and similarly for the bounded variables. 
    
    The multi-set of \emph{free variables} of a derivation $\mathcal{D}$ is denoted by ${\mathsf{fv}}(\mathcal{D})$ and collects all the variables in $\mathcal{D}$ that at some point occur free. Inductively: ${\mathsf{fv}}(A \implies B) := \FV (A)$; 
    ${\mathsf{fv}}\left(\AxiomC{$\mathcal{D}$} \noLine \UnaryInfC{$A\implies B$} \DisplayProof\right) := {\mathsf{fv}}(\mathcal{D}) \cup \FV(A)$ for unary rules and, the free variables of a proof ending on $\wedge \mathcal R$ is defined to be the union of free variables of the two sub-derivations of the antecedents.
\end{definition}

With this standard substitution machinery at hand, we can formulate conditions under which the cut rule is admissible. To do this, we first see that one can define the notion of substitution at the level of proofs.

\begin{definition}
    We define $(A\implies B)\subst{x}{t} \ := \  A\subst{x}{t} \implies B\subst{x}{t}$. For a \qcalc derivation $\mathcal D$ and $y\notin \Big( {\mathsf{ev}}(\mathcal{D}) \cup {\mathsf{bv}}(\mathcal{D}) \Big)$ and for $t$ such that $\{ t\} \cap \Big( {\mathsf{ev}}(\mathcal{D}) \cup {\mathsf{bv}}(\mathcal{D}) \Big) = \varnothing$, the substitution of $t$ for $y$ on $\mathcal D$ is denoted as $\mathcal D \subst{y}{t}$ and is defined recursively as follows 
    \[
    \left(\AxiomC{$\mathcal{D}$}\noLine\UnaryInfC{$A{\implies} B$}\DisplayProof\right)\subst{x}{t}:=\AxiomC{$\mathcal{D}\subst{x}{t}$}\noLine\UnaryInfC{$(A{\implies} B)\subst{x}{t}$}\DisplayProof  \mbox{ and }   
    \]
    \[
    \left(\AxiomC{$\mathcal{D}_0$}
  \AxiomC{$\mathcal{D}_1$}
  \BinaryInfC{$A{\implies} B\land C$}
\DisplayProof\right)\subst{x}{t} :=
\AxiomC{$\mathcal{D}_0\subst{x}{t}$}
  \AxiomC{$\mathcal{D}_1\subst{x}{t}$}
  \BinaryInfC{$(A{\implies} B\land C)\subst{x}{t}$}
\DisplayProof.
 \]   
\end{definition}
It is immediate that the height of a derivation does not change after substitution. The conditions in the  substitution definition warrant that substitutions on proofs yield valid proofs again.

\begin{lemma}\label{theorem:substitutionsOnProofsYieldProofs}
For a \qcalc derivation $\mathcal D$ and $y\notin \Big( {\mathsf{ev}}(\mathcal{D}) \cup {\mathsf{bv}}(\mathcal{D}) \Big)$ and for $t'$ such that $\{ t'\} \cap \Big( {\mathsf{ev}}(\mathcal{D}) \cup {\mathsf{bv}}(\mathcal{D}) \Big) = \varnothing$, the substitution $\mathcal D \subst{y}{t'}$ is again a \qcalc derivation.
\end{lemma}

\begin{proof}
    The proof proceeds by a simple induction of $\mathcal D$. The only non-trivial cases are the two quantifier rules. 
    
    In case the derivation ends as 
    \AxiomC{$\mathcal D$}
    \noLine
    \UnaryInfC{$A\subst{x}{t} \implies B$}
    \RightLabel{\scriptsize $\forall\mathcal L$}
    \UnaryInfC{$\forall x A \implies B$}
    \DisplayProof with $t$ free for $x$ in $A$, the inductive hypothesis warrants that
    \begin{prooftree}
    \AxiomC{$\mathcal D \subst{y}{t'}$}
    \noLine
    \UnaryInfC{$(A\subst{x}{t})\subst{y}{t'} \implies B\subst{y}{t'}$}    
    \end{prooftree}
    Indeed we have access to the IH, since $\Big( {\mathsf{ev}}(\mathcal{D'}) \cup {\mathsf{bv}}(\mathcal{D'}) \Big)\subseteq \Big( {\mathsf{ev}}(\mathcal{D}) \cup {\mathsf{bv}}(\mathcal{D}) \Big)$ for $\mathcal{D'}$ a sub-proof of $\mathcal{D}$. 
    
    Since $\big(A\subst{x}{t}\big)\subst{y}{t'} = \big(A\subst{y}{t'}\big)\subst{x}{\big(t\subst{y}{t'}\big)}$ by Proposition \ref{theorem:NestedSubstitution} and by the assumptions on $y$ and $t'$ we see that $t\subst{y}{t'}$ is free for $x$ in $A\subst{y}{t'}$ and $x\neq y$. The same assumptions on $y$ and $t'$ yield that $t'$ is free for $y$ in $B$. Consequently, we can apply the {$\forall\mathcal L$} rule again to conclude $\forall x A \subst{y}{t'} \implies B \subst{y}{t'}$ as was to be shown. 

    The {$\forall\mathcal R$} rule goes similarly but is easier.
\end{proof}

With this definition and lemma we can prove a controlled form of cut-admissibility for \qcalc.

\begin{theorem}[Cut admissibility for \qcalc]
\label{thm:qcalc-cut}
If
$\mathcal{D_0}\vdash_\qcalc A\implies B$ and
$\mathcal{D_1}\vdash_\qcalc B\implies C$, so that ${\mathsf{fv}}(\mathcal{D_1}) \cap \big( {\mathsf{ev}}(\mathcal{D_0}) \cup {\mathsf{bv}}(\mathcal{D_0})\big) = \varnothing $ then, there is some $\mathcal{D_2}$ with 
$\mathcal{D_2}\vdash_{\qcalc} A\implies C$ with ${\mathsf{ev}}(\mathcal{D_2})\subseteq {\mathsf{ev}}(\mathcal{D_0}) \cup {\mathsf{ev}}(\mathcal{D_1})$ with ${\mathsf{bv}}(\mathcal{D_2})\subseteq {\mathsf{bv}}(\mathcal{D_0}) \cup {\mathsf{bv}}(\mathcal{D_1})$ and with ${\mathsf{fv}}(\mathcal{D_2})\subseteq {\mathsf{fv}}(\mathcal{D_0}) \cup {\mathsf{fv}}(\mathcal{D_1})$.
\end{theorem}

\begin{proof}
    We repeat the proof of Theorem \ref{thm:calc-cut} keeping track of the variable conditions and only have to pay special attention to the case that $B=\forall x B_0$ and principal on both sides (disjointness conditions are preserved under taking subderivations). 
    
    The premisses to $A\implies \forall x B_0$ and $\forall x B_0 \implies C$ are 
    \begin{center}
        \AxiomC{$\mathcal{D}_{00}$}\noLine \UnaryInfC{$A\implies B_0\subst{x}{y}$}\DisplayProof and
    \AxiomC{$\mathcal{D}_{10}$}\noLine \UnaryInfC{$ B_0\subst{x}{t} \implies C$}\DisplayProof    
    \end{center}
     respectively for a certain term $t$ and variable $y$. Since $t$ has empty intersection with both the eigenvariables of $\mathcal{D}_{00}$ and the bound variables of $\mathcal{D}_{00}$ we can substitute $t$ for $y$ in $\mathcal{D}_{00}$ which yields (by Lemma \ref{theorem:substitutionsOnProofsYieldProofs} and since $y\notin A\wedge \forall x B_0$) a proof  of $A\implies B_0\subst{x}{t}$. Since the condition on free versus bound and eigenvariables holds again at the level of $\mathcal{D}_{00}$ and  $\mathcal{D}_{10}$, we can apply the inductive hypothesis.
\end{proof}

As in the standard treatments of cut-elimination, we have seen that in \textit{sensu stricto} we do not have cut admissibility but if the bound variables and eigenvariables are chosen suitably, we do. The easiest way around the possible variable clashes is by hardwiring it in the syntax and requiring that bound and free variables come from different alphabets. In most of the literature however, one takes the notion of $\alpha$-conversion to be central as we shall define below. Basically, $\alpha$-equivalent formulas are equivalent for all essential purposes and only differ in the usage of bound variables.

\begin{definition}
    The $\alpha$-equivalence relation $\equiv_\alpha$ is the smallest equivalence relation on $\qlang$ formulas such that if $y$ is free for $x$ in $A$, and if $A':= A \subst{x}{y}$ and if $\FV(\forall x A) = \FV (\forall y A')$, then $\forall x A \equiv_\alpha \forall y A'$. Moreover, $A\equiv_\alpha A'$ and  $B\equiv_\alpha B'$ implies $A\wedge B\equiv_\alpha A'\wedge B'$ and, $\Diamond A\equiv_\alpha \Diamond A'$.
\end{definition}

We recall that in this paper, $\FV(A)$ is a multiset and not just a set. Consequently, we see, for example, $\forall x R(x,y) \wedge S(y) \not \equiv_\alpha \forall y R(y,y) \wedge S(y)$, as is to be expected. It is easy to see that each formula $A$ is $\alpha$-equivalent to one where the bound variables are chosen from a (large enough) set $V$ of variables such that each bound variable is bound by exactly one quantifier, that is, has multiplicity one in the multi-set. The same holds for proofs.\footnote{Here we are making explicit the machinery left implicit in the cut elimination arguments proposed by, e.g.,~\cite{NegriVonPlato:2001:StructuralProofTheory} disallowing formula occurrences like $\forall x P(x) \wedge \forall x Q(x)$.}

\begin{lemma}\label{theorem:pureVariableProofs}
    For each derivation \AxiomC{$\mathcal{D}$}\noLine \UnaryInfC{$A\implies B$}\DisplayProof and, for each $A\equiv_\alpha A'$ and $B\equiv_\alpha B'$ so that all bound variables in $A'\wedge B'$ occur with multiplicity 1, for each large enough disjoint sets $V$ and $W$ of variables, we can find a derivation \AxiomC{$\mathcal{D'}$}\noLine \UnaryInfC{$A'\implies B'$}\DisplayProof of equal height where all bound and eigenvariables come from $V$ and $W$ respectively and all occur with multiplicity 1.
\end{lemma}

\begin{proof}
     The proof proceeds by induction on $\mathcal{D}$. 
\end{proof}

We say that a derivation is a \textit{pure-variable derivation}  whenever the free variables, the eigenvariables and the bound variables are disjoint and where the latter two only occur with multiplicity one as in Lemma \ref{theorem:pureVariableProofs}. Pure variable derivations do allow cut-elimination.

\begin{theorem}\label{theorem:CutAdmissibilityQRCOne}
    Any pure-variable derivation $\AxiomC{$\mathcal{D}$}\noLine\UnaryInfC{$A\implies B$}\DisplayProof$ of $\qcalc + \textit{Cut}$ can be effectively transformed into a \qcalc derivation $\AxiomC{$\mathcal{D'}$}\noLine\UnaryInfC{$A\implies B$}\DisplayProof$.
\end{theorem}

\begin{proof}
    The procedure consists of removing cuts one-by-one top-down, using Lemma \ref{thm:qcalc-cut}. The pure-variable condition together with the variable conditions of Lemma \ref{thm:qcalc-cut} ensure that we can repeatedly apply the Lemma.
\end{proof}

\subsection{Completeness}

In Theorem \ref{thm:modular-soundness} we have seen that every theorem of \calc is also a theorem of $\RC_1$. We shall now see that the converse also hold so that both calculi prove the exact same set of theorems.

\begin{theorem}\label{theorem:CalcEquivToRCOne}
    For \lang formulas $A$ and $B$ we have
    $\calc \vdash A\implies B$ if and only if $A\vdash_{\RC_1}B$.
\end{theorem}

\begin{proof}
    In virtue of Theorem \ref{thm:modular-soundness} we only need to show that if $A\vdash_{\RC_1}B$, then we have $\calc \vdash A\implies B$. The proof follows by induction on an $\RC_1$ derivation. 
    
    We first focus on the axioms. 
    
    The $  \AxiomC{}
  \RightLabel{\scriptsize $\top$}
  \UnaryInfC{$A\vdash_{\RC_1}\top$}
  \DisplayProof
$ axiom is present in \calc and the $  \AxiomC{}
  \RightLabel{\scriptsize id}
  \UnaryInfC{$A\vdash_{\RC_1}A$}
  \DisplayProof$ axiom is fine since we proved by induction on $A$ that $\calc \vdash A\implies A$ (Lemma \ref{lem:modular-generalised-id}). 
  
  The $\AxiomC{}
  \RightLabel{\scriptsize $\land\mathcal L_1$}
  \UnaryInfC{$A\land B\vdash_{\RC_1}A$}\DisplayProof$ 
  axiom of $\RC_1$ can be obtained by the \calc reasoning starting from a \calc derivation of $A\implies A$ followed by 
  $\AxiomC{$A\implies A$}\RightLabel{\scriptsize $\land\mathcal L_1$}\UnaryInfC{$A\wedge B\implies A$}\DisplayProof$ and similarly for {\scriptsize $\land\mathcal L_2$}. 
  
    The $    \AxiomC{}
    \RightLabel{\scriptsize 4}
    \UnaryInfC{$\Diamond\Diamond A\vdash_{\RC_1}\Diamond A$}
    \DisplayProof$
    axiom in $\RC_1$ can be obtained by the \calc reasoning starting from a \calc derivation of $\Diamond A\implies \Diamond A$ followed by $\AxiomC{$\Diamond A\implies \Diamond A$}    \RightLabel{\scriptsize $4\mathrm{Nec}$}
    \UnaryInfC{$\Diamond \Diamond A\implies\Diamond B$}
    \DisplayProof$.

    For the rules, we observe that the $\RC_1$ rules {\scriptsize Nec} and {\scriptsize $\land\mathcal R$} are simply present as rules in translated form in \calc. Finally, the {\scriptsize Cut} rule is covered by cut-admissibility of \calc derivations (Theorem \ref{thm:calc-cut}).
\end{proof}

As already observed above, $\qcalc \nvdash \forall x\, \forall y \ \big(P(x) \wedge Q(y) \big) \implies P(y)$ yet it easy to see that $\forall x\, \forall y \ \big(P(x) \wedge Q(y) \big) \vdash_{\QRC} P(y)$ so, we cannot expect full proof-equivalence between $\QRC$ and \qcalc. If we consider theorems modulo $\alpha$-equivalence, then indeed we do have proof equivalence between the two calculi. We first need a result to the extent that $\QRC$ derivations also behave well w.r.t.~$\alpha$-conversion.

\begin{lemma}\label{theorem:QRCprovabilityAlphaInvariant}
    If $A\vdash_{\QRC}B$ and if $A\equiv_\alpha A'$ and if $B\equiv_\alpha B'$, then $A'\vdash_{\QRC}B'$. Moreover, both derivations have equal heights.
\end{lemma}

\begin{proof}
    By induction on the derivation of $A\vdash_{\QRC}B$.
\end{proof}

\begin{theorem}\label{theorem:QRC1andQlangEquivalent}
    Let $A,B$ be $\qlang$ formulas with all bounded variables occuring with multiplicity one, and where $ \BV (A\wedge B) \cap \FV(A\wedge B) = \varnothing$.  Under these conditions, we have $\qcalc \vdash A\implies B$ if and only if $A\vdash_{\QRC}B$.    
\end{theorem}

\begin{proof}
    The direction, if $\qcalc\vdash A\implies B$, then $A\vdash_{\QRC} B$, does not depend on any assumptions on the variables in $A$ and $B$ and has been proven in Theorem \ref{thm:modular-soundness}. 
    
    The other direction is proven by induction on a $\QRC$ derivation and we only have to reconsider the {Cut} rule beside considering the quantifier rules. 
    
    The reasoning for the other axioms and rules lift along the proof of Theorem \ref{theorem:CalcEquivToRCOne}.

    We first consider the $\AxiomC{$A\vdash_{\QRC}B$}
  \AxiomC{$B\vdash_{\QRC}C$}
  \RightLabel{\scriptsize Cut}
  \BinaryInfC{$A\vdash_{\QRC}C$}
  \DisplayProof$-rule. Choose $B'\equiv_\alpha B$, so that we have access to the inductive hypothesis. In particular, we choose $B'$ so that
  \[
BV(A)\cap\BV(B') = BV(B')\cap\BV(C) = \varnothing = \BV (A\wedge B') \cap \FV(A\wedge B') = \BV (B'\wedge C) \cap \FV(B'\wedge C)
  \]
and, the multiplicity of each bound variable in $A\wedge B'\wedge C$ is one. Renaming bounded variables does not change the $\QRC$ provability status (Lemma \ref{theorem:QRCprovabilityAlphaInvariant}) so that by the  inductive hypotheses we obtain $\mathcal{D_0}\vdash_{\qcalc} A \implies B$ and $\mathcal{D_1}\vdash_{\qcalc} B \implies C$. The eigenvariables in $\mathcal{D}_i$ can be chosen (using Lemma \ref{theorem:pureVariableProofs} taking $A'=A$ and $B'=B$) so that the conditions for a single cut admissibility (Theorem \ref{thm:qcalc-cut}) are satisfied. Consequently, we obtain $\qcalc \vdash A\implies C$ as was to be shown.

We now consider the quantifier rules. The {\scriptsize $\forall\mathcal L$} rule in \QRC is the same in translated form in \qcalc. Likewise, the \AxiomC{$A\vdash_{\QRC}B$}
  \RightLabel{\scriptsize $\forall\mathcal R$}
  \UnaryInfC{$A\vdash_{\QRC}\forall x\,B$}
  \DisplayProof with {\mbox{{ $x\notin\FV(A)$}}} is covered by 
  \AxiomC{$A\implies B\subst{x}{y}$}
  \RightLabel{\scriptsize $\forall\mathcal R$}
  \UnaryInfC{$A\implies\forall x\,B$}
  \DisplayProof with $y$ free for $x$ in $B$; $y\notin\FV(A)\cup\FV(\forall x\,B)$
since, clearly, $x$ is free for $x$ in $B$ and, $x\notin \FV (\forall x B)$. 

The rules
$\AxiomC{$A\vdash_{\QRC}B$}
  \RightLabel{\scriptsize Inst}
  \UnaryInfC{$A\subst{x}{t}\vdash_{\QRC}B\subst{x}{t}$}
  \DisplayProof$
({ $t$ free for $x$ in both $A$ and $B$}), and 
  $\AxiomC{$A\subst{x}{c}\vdash_{\QRC}B\subst{x}{c}$}
  \RightLabel{\scriptsize CE}
  \UnaryInfC{$A\vdash_{\QRC}B$}
  \DisplayProof$
  where {$c$ occurs in neither $A$ nor $B$} are easily seen to be height-preserving admissible in \qcalc.\footnote{That is, 
  by an easy induction on \qcalc derivations we see that if
  $\qcalc \vdash_{\leq n} A \implies B$, then \mbox{$\qcalc \vdash_{\leq n} A\subst{x}{t}\implies B\subst{x}{t}$} whenever $t$ free for $x$ in both $A$ and $B$. Similarly, we can prove that if
  $\qcalc \vdash_{\leq n} A\subst{x}{c}\implies B\subst{x}{c}$, then $\qcalc \vdash_{\leq n} A\implies B$ whenever 
    $c$ occurs in neither $A$ nor $B$.}
    \end{proof}

\section{Cut elimination in the non-analytic calculi}
\label{sec:cut-non-analytic}

In the previous section, we have proven that both \calc and \qcalc satisfy a certain form of cut admissibility. In Theorem \ref{thm:modular-gentzen-equivalence} 
we have shown that \calc,
\calcFourNec, and \calcFourNecStar all prove the same set of theorems and likewise for the quantified versions. As a consequence, we get cut-admissibility for those systems as well. 

In this section we will showcase how the strict-positiveness makes up for non-analyticity of the proof systems \calcFourNec by giving a direct proof of cut-elimination \calcFourNec. We start with some preparatory results.

\begin{lemma}[Top replacement]
\label{lem:g1k4-topNew}
If
$\mathsf \calcFourNec \vdash_{\leq n}\top\implies C$, then, for every $A$ in \lang, $\mathsf \calcFourNec \vdash_{\leq n}A\implies C$.
\end{lemma}

\begin{proof}
By an easy induction on a derivation of $\top\implies C$.  The possible final rules are $\top\mathcal R$, $\land\mathcal R$,
$4\mathcal R$. 
In all cases, we apply the induction hypothesis to the premise or
premises and reconstruct the final rule. 
\end{proof}

The next lemma is standard in structural proof theory.

\begin{lemma}[Height-preserving invertibility of $\land\mathcal R$]
\label{lem:and-r-inv}
If $\calcFourNec\vdash_{\leq n} A\implies B\land C$,
then
$\calcFourNec\vdash_{\leq n} A\implies B$ and
$\calcFourNec\vdash_{\leq n} A\implies C$.
\end{lemma}

\begin{proof}
By induction on the height of the derivation. If its last rule is $\land\mathcal R$, the premises give the result. If the last rule is $\land\mathcal L_i$ or $4\mathcal L$, apply the induction hypothesis to the premise and reapply the same rule. No other rule can have a conjunction as its consequent.
\end{proof}

The notion of \textit{principal formula} in \calcFourNec is as in \calc with the further stipulations that  
$\Diamond\Diamond A$ is principal in $4\mathcal L$, and the
consequent $\Diamond B$ is principal in $4\mathcal R$.  

\begin{lemma}[Modal conjunction extension]
\label{lem:g1k4-diamond-conjunctionNew}
The following
two rules are admissible in $\calcFourNec$:
\begin{center}
\begin{tabular}{c@{\qquad\qquad}c}
  \AxiomC{$\Diamond A\implies C$}
  \UnaryInfC{$\Diamond(A\land B)\implies C$}
  \DisplayProof
&
  \AxiomC{$\Diamond A\implies C$}
  \UnaryInfC{$\Diamond(B\land A)\implies C$}
  \DisplayProof
\end{tabular}
\end{center}
\end{lemma}

\begin{proof}
It suffices to prove the first rule; the second is symmetric.  We let
$\mathcal{D}\vdash\Diamond A\implies C$ and use lexicographic induction on
$\langle h(\mathcal{D}),|C|\rangle$, with derivation height primary.

If $\mathcal{D}$ ends in $\top\mathcal R$, the claim is immediate.  

If it ends in Nec, then $C=\Diamond C_0$ and its premise is $A\implies C_0$;
we apply $\land\mathcal L_1$ to the premise and then Nec.  

If it ends in
$4\mathcal R$ or $\land\mathcal R$, we apply the primary induction
hypothesis to the premise or premises and reconstruct the same right rule.

It remains to consider a final $4\mathcal L$.  This implies that $A=\Diamond A_0$
and the premise is $\Diamond A_0\implies C$.  

If $C=\top$, the claim is immediate.  

If $C=\Diamond C_0$, we use
\begin{prooftree}
  \AxiomC{$\Diamond A_0\implies\Diamond C_0$}
  \RightLabel{\scriptsize $\land\mathcal L_1$}
  \UnaryInfC{$\Diamond A_0\land B\implies\Diamond C_0$}
  \RightLabel{\scriptsize Nec}
  \UnaryInfC{$\Diamond(\Diamond A_0\land B)
    \implies\Diamond\Diamond C_0$}
  \RightLabel{\scriptsize $4\mathcal R$}
  \UnaryInfC{$\Diamond(\Diamond A_0\land B)\implies\Diamond C_0$}
\end{prooftree}

If $C=C_1\land C_2$, Lemma~\ref{lem:and-r-inv} gives
$\Diamond A_0\implies C_i$ with no increase in height.  We reapply
$4\mathcal L$ to obtain $\Diamond\Diamond A_0\implies C_i$ and then use
the secondary induction hypothesis, since $|C_i|<|C|$, to obtain
$\Diamond(\Diamond A_0\land B)\implies C_i$.  We finish with
$\land\mathcal R$.

The remaining possibility that $C$ is atomic cannot occur: with a modal antecedent,
no rule producing an atomic consequent is available except repeated
$4\mathcal L$, and following those premises upwards would eventually
require a zero-premise derivation with a modal antecedent and atomic
consequent.
\end{proof}

We are now ready to prove Cut admissibility.

\begin{theorem}[Cut admissibility for \calcFourNec]
\label{thm:g1k4-cut}
$\calcFourNec\vdash A\implies C$ and
$\calcFourNec\vdash C\implies B$, then
$\calcFourNec\vdash A\implies B$.
\end{theorem}

\begin{proof} We give a constructive proof of admissibility of the cut rule.

Let cut-free derivations
$\mathcal{D}_l\vdash A\implies C$ and $\mathcal{D}_r\vdash C\implies B$ be given, and put
$h_l=h(\mathcal{D}_l)$ and $h_r=h(\mathcal{D}_r)$.  We argue by lexicographic induction
on $\langle h_l,h_r\rangle$, with the left component primary: no complexity measure on the cut formula is required.

If one premise has height zero, the Initial cases are immediate and the
$\top\mathcal R$ cases are immediate except when the left derivation
concludes $A\implies\top$; then
Lemma~\ref{lem:g1k4-topNew} transforms the right derivation
$\top\implies B$ into one of $A\implies B$.

Suppose next that the last inference of $\mathcal{D}_l$ leaves the cut formula in
its succedent unchanged:  it is then $\land\mathcal L_i$ or $4\mathcal L$.
Thus, we can push the cut above that inference, apply the induction hypothesis, and
apply the same left rule.  The left height strictly decreases.

Dually, if the last inference of $\mathcal{D}_l$ introduces the outer connective
of $C$ but the last inference of $\mathcal{D}_r$ leaves its antecedent unchanged,
we push the cut above the final right inference.  The possibilities are
$\land\mathcal R$ and $4\mathcal R$; we may use one call for
$4\mathcal R$ and two for $\land\mathcal R$.  The left height is fixed
and each right height decreases.

It remains to treat cuts for which $C$ is principal in both final
inferences.

\begin{itemize}
\item If $C=C_1\land C_2$, the final rules are
$\land\mathcal R/\land\mathcal L_i$.  Then, we cut the $i$-th left premise
$A\implies C_i$ against the premise $C_i\implies B$ of the right rule.
The left height strictly decreases.

\item If the final rules are Nec/Nec, we replace the cut on
$\Diamond C_0$ by the cut between their premises on $C_0$, and then
reapply Nec.

\item For $4\mathcal R/4\mathcal L$, the configuration is
\begin{prooftree}
  \AxiomC{$A\implies\Diamond\Diamond\Diamond C_0$}
  \RightLabel{\scriptsize $4\mathcal R$}
  \UnaryInfC{$A\implies\Diamond\Diamond C_0$}
  \AxiomC{$\Diamond C_0\implies B$}
  \RightLabel{\scriptsize $4\mathcal L$}
  \UnaryInfC{$\Diamond\Diamond C_0\implies B$}
  \RightLabel{\scriptsize Cut}
  \BinaryInfC{$A\implies B$}
\end{prooftree}
Then we apply $4\mathcal L$ to the right premise, obtaining
$\Diamond\Diamond\Diamond C_0\implies B$, and we cut this against the
premise of the left $4\mathcal R$.  The new pair has measure
$\langle h_l-1,h_r+1\rangle$.

\item For $4\mathcal R$/Nec, the configuration is
\begin{prooftree}
  \AxiomC{$A\implies\Diamond\Diamond C_0$}
  \RightLabel{\scriptsize $4\mathcal R$}
  \UnaryInfC{$A\implies\Diamond C_0$}
  \AxiomC{$C_0\implies B_0$}
  \RightLabel{\scriptsize Nec}
  \UnaryInfC{$\Diamond C_0\implies\Diamond B_0$}
  \RightLabel{\scriptsize Cut}
  \BinaryInfC{$A\implies\Diamond B_0$}
\end{prooftree}
Then, we apply $4\mathcal L$ to the right conclusion, obtaining
$\Diamond\Diamond C_0\implies\Diamond B_0$, and we cut this against the
premise of the left $4\mathcal R$.  The new measure is
$\langle h_l-1,h_r+1\rangle$.

\item The remaining principal modal pair is Nec/$4\mathcal L$:
\begin{prooftree}
  \AxiomC{$A_0\implies\Diamond C_0$}
  \RightLabel{\scriptsize Nec}
  \UnaryInfC{$\Diamond A_0\implies\Diamond\Diamond C_0$}
  \AxiomC{$\Diamond C_0\implies B$}
  \RightLabel{\scriptsize $4\mathcal L$}
  \UnaryInfC{$\Diamond\Diamond C_0\implies B$}
  \RightLabel{\scriptsize Cut}
  \BinaryInfC{$\Diamond A_0\implies B$}
\end{prooftree}
Let $\mathcal{P}\vdash A_0\implies\Diamond C_0$ be the derivation of the premise of the final Nec.
Its height is $h_l-1$.  Since its consequent is modal, $\mathcal{P}$ is not a
zero-premise derivation and cannot end in $\land\mathcal R$.  We need to inspect
its possible last rule.

If $\mathcal{P}$ ends in $4\mathcal R$, its premise is
$A_0\implies\Diamond\Diamond C_0$.  Then we apply Nec to that premise, obtaining
$\Diamond A_0\implies\Diamond\Diamond\Diamond C_0$ with height
$h_l-1$; we apply $4\mathcal L$ twice to the right premise
$\Diamond C_0\implies B$, and then use the induction hypothesis.  The left
height has decreased.

If $\mathcal{P}$ ends in Nec, then $A_0=\Diamond A_1$ and its premise is
$A_1\implies C_0$.  That is: $\mathcal{P}$ itself is
$\Diamond A_1\implies\Diamond C_0$. We cut $\mathcal{P}$ against the premise
$\Diamond C_0\implies B$ of the final right $4\mathcal L$; the left
height is $h_l-1$ and the right height is $h_r-1$. Then we reapply
$4\mathcal L$ to obtain $\Diamond\Diamond A_1\implies B$.

If $\mathcal{P}$ ends in $\land\mathcal L_i$, let us write
$A_0=A_1\land A_2$ and let its premise be
$A_i\implies\Diamond C_0$.  We apply Nec to obtain
$\Diamond A_i\implies\Diamond\Diamond C_0$, cut this against the full
right derivation $\Diamond\Diamond C_0\implies B$, and use
Lemma~\ref{lem:g1k4-diamond-conjunctionNew} to pass from
$\Diamond A_i\implies B$ to $\Diamond(A_1\land A_2)\implies B$.
The new cut has left height $h_l-1$.

Finally, if $\mathcal{P}$ ends in $4\mathcal L$, then
$A_0=\Diamond\Diamond A_1$ and its premise is
$\Diamond A_1\implies\Diamond C_0$.  We cut that premise against
$\Diamond C_0\implies B$ and then apply $4\mathcal L$ twice.  Again the
left height strictly decreases.
\end{itemize}

These cases exhaust all pairs of final rules.
\end{proof}

The following corollary follows easily.

\begin{corollary}[Effective cut elimination in \calcFourNec]
\label{cor:g1k4-cut-elimination}
Adding Cut to \calcFourNec does not enlarge its derivability relation.
Every derivation in the extended calculus can be effectively transformed
into a cut-free \calcFourNec-derivation.
\end{corollary}

\begin{proof}
Recursively eliminate applications of cut in the immediate subderivations and then use
Theorem~\ref{thm:g1k4-cut} on the resulting cut-free premises.
\end{proof}

The theorem can easily be extended to the quantified setting, but we only included the results of this section as an example of how poor calculi can still work in the context of strictly positive fragments.

\section{Conservativity over \texorpdfstring{$\gl$}{GL}}
\label{sec:gl-positive-fragment}

Dashkov proved that $\RC_1$ is exactly the strictly positive fragment of
$\gl$~\cite{Dashkov:2012:PositiveFragment}.  
In this section, we provide a novel proof of that conservativity result by a purely syntactic extraction argument.  

We fix first from the outset a standard G3-style multiset calculus ($\SGL$) whose
primitive propositional language contains $\neg$, $\land$, and $\top$.
Cut elimination for ordinary multiset sequent presentations of $\gl$
is established in~\cite{GoreRamanayake:2012:Valentini}; Lyon's recent
analysis makes explicit the relation between such sequent systems and
several other structural formalisms for $\gl$~\cite{Lyon:2025:UnifyingGL}.\footnote{Direct proof systems for richer positive fragments of $\gl$
have recently been studied in~\cite{Tanaka:2026:PositiveGL}.  For
comparison, cut-free sequent analyses are also available for the
intuitionistic G\"odel--L\"ob logic and for intuitionistic strong L\"ob
logic~\cite{VanDerGiessenIemhoff:2021:SequentIGL,
ShillitoEtAl:2023:StrongLob}; these systems are not used in the extraction
argument below, but they locate the multiset treatment within the current
proof-theoretic landscape around L\"ob principles.}
Our point here is to show how a cut-free multiset derivation in $\gl$ can
be collapsed to a derivation in our sequent calculus
\calc, over a shared grammar.

Let $\mathcal L_{\Box}$ be the classical modal language generated by
propositional variables, $\top$, $\neg$, $\land$, and $\Box$. The multiset calculus \SGL is defined by the rules displayed in
\Cref{fig:sgl-rules}.  It is the usual two-sided G3 presentation of $\gl$ in
the Boolean basis $\{\neg,\land,\top\}$, with weakening absorbed into the
initial and modal rules.  Lyon works explicitly with finite multisets in his
system $\mathsf{GL}_{\mathrm{seq}}$, records soundness, completeness,
syntactic cut elimination, and admissibility of weakening and contraction,
and explains how that system is related by proof transformations to the
cyclic, non-wellfounded, tree-hypersequent, labelled, and linear-nested
systems for $\gl$~\cite[Sec.~5.1]{Lyon:2025:UnifyingGL}.  His chosen Boolean
basis is $\{\neg,\lor,\Box\}$; the basis used here is convenient for a direct comparison with our system ~\calc based on \lang.

\begin{figure}[H]
\centering
\begin{tabular}{c@{\qquad}c}
  \AxiomC{}
  \RightLabel{\scriptsize id}
  \UnaryInfC{$\Gamma,p\Rightarrow p,\Delta$}
  \DisplayProof
&
  \AxiomC{}
  \RightLabel{\scriptsize $\top\mathcal R$}
  \UnaryInfC{$\Gamma\Rightarrow\top,\Delta$}
  \DisplayProof
\\[3ex]
  \AxiomC{$\Gamma\Rightarrow A,\Delta$}
  \RightLabel{\scriptsize $\neg\mathcal L$}
  \UnaryInfC{$\neg A,\Gamma\Rightarrow\Delta$}
  \DisplayProof
&
  \AxiomC{$A,\Gamma\Rightarrow\Delta$}
  \RightLabel{\scriptsize $\neg\mathcal R$}
  \UnaryInfC{$\Gamma\Rightarrow\neg A,\Delta$}
  \DisplayProof
\\[3ex]
  \AxiomC{$A,B,\Gamma\Rightarrow\Delta$}
  \RightLabel{\scriptsize $\land\mathcal L$}
  \UnaryInfC{$A\land B,\Gamma\Rightarrow\Delta$}
  \DisplayProof
&
  \AxiomC{$\Gamma\Rightarrow A,\Delta$}
  \AxiomC{$\Gamma\Rightarrow B,\Delta$}
  \RightLabel{\scriptsize $\land\mathcal R$}
  \BinaryInfC{$\Gamma\Rightarrow A\land B,\Delta$}
  \DisplayProof
\end{tabular}
\medskip
\begin{prooftree}
  \AxiomC{$\Box\Gamma,\Gamma,\Box C\Rightarrow C$}
  \RightLabel{\scriptsize $\Box$-L\"ob}
  \UnaryInfC{$\Omega,\Box\Gamma\Rightarrow\Box C,\Delta$}
\end{prooftree}
\caption{The multiset calculus \SGL.  The occurrence of $\Box C$ is called the diagonal formula of the modal rule.}
\label{fig:sgl-rules}
\end{figure}

Taking $\Diamond$ as primitive via the classical modal duality
 $\Diamond A\equiv\neg\Box\neg A$,
the corresponding full $\Diamond$-presentation, denoted
$\SGLDiamond_{\text{full}}$, has the same classical G3 rules of \SGL
and replaces $\Box$-L\"ob by
\begin{prooftree}
  \AxiomC{$A\Rightarrow\Sigma,\Diamond\Sigma,\Diamond A$}
  \RightLabel{\scriptsize $\Diamond$-L\"ob}
  \UnaryInfC{$\Omega,\Diamond A\Rightarrow\Diamond\Sigma,\Delta$}
\end{prooftree}
where $\Diamond\Sigma$ is the occurrencewise image of the multiset $\Sigma$.

We now restrict this presentation to the strictly positive language.

\begin{definition}[The positive multiset calculus $\SGLDiamond$]
\label{def:sgl-diamond}
The calculus $\SGLDiamond$ is the restriction of
$\mathsf{SGL}^{\Diamond}_{\mathrm{full}}$ to sequents over \lang.  Its rules
are those in \Cref{fig:sgl-diamond-rules}.
\end{definition}

\begin{figure}[H]
\centering
\begin{tabular}{c@{\qquad\qquad}c}
  \AxiomC{}
  \RightLabel{\scriptsize id}
  \UnaryInfC{$\Gamma,p\Rightarrow p,\Delta$}
  \DisplayProof
&
  \AxiomC{}
  \RightLabel{\scriptsize $\top\mathcal R$}
  \UnaryInfC{$\Gamma\Rightarrow\top,\Delta$}
  \DisplayProof
\\[3ex]
  \AxiomC{$A,B,\Gamma\Rightarrow\Delta$}
  \RightLabel{\scriptsize $\land\mathcal L$}
  \UnaryInfC{$A\land B,\Gamma\Rightarrow\Delta$}
  \DisplayProof
&
  \AxiomC{$\Gamma\Rightarrow\Delta,A$}
  \AxiomC{$\Gamma\Rightarrow\Delta,B$}
  \RightLabel{\scriptsize $\land\mathcal R$}
  \BinaryInfC{$\Gamma\Rightarrow\Delta,A\land B$}
  \DisplayProof
\end{tabular}
\medskip
\begin{prooftree}
  \AxiomC{$A\Rightarrow\Sigma,\Diamond\Sigma,\Diamond A$}
  \RightLabel{\scriptsize $\Diamond$-L\"ob}
  \UnaryInfC{$\Omega,\Diamond A\Rightarrow\Diamond\Sigma,\Delta$}
\end{prooftree}
\caption{The strictly positive fragment $\SGLDiamond$ for \gl The occurrence of $\Diamond A$ is called the diagonal formula of the modal rule.}
\label{fig:sgl-diamond-rules}
\end{figure}

\begin{lemma}[Positive restriction]\label{lem:gl-positive-restriction}
    Let $\Gamma,\Delta$ be finite multisets of formulas in \lang.  If $\SGLDiamond_{\text{full}}\vdash\Gamma\Rightarrow\Delta$, then $\SGLDiamond\vdash\Gamma\Rightarrow\Delta$.
\end{lemma}
\begin{proof}
    Every rule of $\SGLDiamond_{\text{full}}$ is analytic. Thus, every formula in the derivation is strictly positive: no negation rule occurs, and every inference is one of the rules of $\SGLDiamond$.
\end{proof}

\begin{corollary}[Positive cut-free completeness]\label{cor:gl-positive-cut-free-completeness}
    For $A,B\in\lang$, if $\gl\vdash A\to B$, then $\SGLDiamond\vdash A\Rightarrow B$.
\end{corollary}

\subsection{Preliminary proof analysis}\label{subsec:gl-positive-extraction}

For a finite positive multiset $\Gamma$, define
\[
 \mdepth(\Gamma):=
 \max\{\mdepth(A):A\text{ occurs in }\Gamma\},
 \qquad
 \mdepth(\varnothing):=0.
\]

\begin{lemma}[Positive modal-depth witness]\label{lem:gl-positive-depth}
    If $\SGLDiamond\vdash\Gamma\Rightarrow\Delta$, then some occurrence $D$ in $\Delta$ satisfies
    $$\mdepth(D)\leq\mdepth(\Gamma).$$
    Consequently, $\SGLDiamond\nvdash A\Rightarrow\Diamond A$.
\end{lemma}
\begin{proof}
    By induction on the derivation height. 
    
    The cases for height zero are immediate. 

    The rule $\land\mathcal L$ preserves the maximum modal depth of the antecedent.

    If the last rule is $\land\mathcal R$, we apply the induction hypothesis to both premises. If either selected occurrence lies in the common side multiset $\Delta$, we keep that occurrence. Otherwise the two selected occurrences are the active formulas $A$ and $B$, and we select $D:=A\land B$; its modal depth is the maximum of their modal depths.

    For the modal rule
        \begin{prooftree}
            \AxiomC{$A\Rightarrow\Sigma,\Diamond\Sigma,\Diamond A$}
            \RightLabel{\scriptsize $\Diamond$-L\"ob}
            \UnaryInfC{$\Omega,\Diamond A\Rightarrow\Diamond\Sigma,\Delta$}
        \end{prooftree}
    the induction hypothesis yields an occurrence $D'$ in the succedent of the premise such that $\mdepth(D')\leq\mdepth(A)$. Hence $D'$ cannot be the distinguished occurrence of $\Diamond A$.  If $D'=S$ is an occurrence in $\Sigma$, we select the corresponding occurrence $\Diamond S$ in the conclusion. Then we have $\mdepth(\Diamond S)\leq \mdepth(\Diamond A)\leq \mdepth(\Omega,\Diamond A)$. If instead $D'=\Diamond S$ is an occurrence in $\Diamond\Sigma$, we keep that occurence and have $\mdepth(\Diamond S)\leq\mdepth(A)<\mdepth(\Diamond A)$.

    The final claim is obtained by direct application of such an invariant to the sequent $A\Rightarrow\Diamond A$.
\end{proof}

\begin{lemma}[Weakening in $\SGLDiamond$]\label{lem:gl-positive-weakening}
    If $\SGLDiamond\vdash\Gamma\Rightarrow\Delta$, then $\SGLDiamond\vdash\Pi,\Gamma\Rightarrow\Delta,\Theta$ for arbitrary finite multisets $\Pi,\Theta$.
\end{lemma}
\begin{proof}
    Straightforward induction on the height of the derivation.
\end{proof}

If $\Gamma$ is a finite multiset, let
$\bigwedge\Gamma$ be a fixed bracketing of the conjunction of all its
formula occurrences, and let us define $\bigwedge\varnothing:=\top$.

The following is a technical lemma.

\begin{lemma}[Conjunctive contexts]\label{lem:gl-conjunctive-contexts}
    In \calc, conjunctions obtained from different enumerations and bracketings
    of the same finite multiset are interderivable.  Moreover:
    \begin{enumerate}[label=\upshape(\roman*)]
      \item if an occurrence of $C$ belongs to $\Gamma$ and
      $\calc\vdash C\implies D$, then
      $\calc\vdash\bigwedge\Gamma\implies D$;
      \item $\bigwedge(A\land B,\Gamma)$ and
      $\bigwedge(A,B,\Gamma)$ are interderivable.
    \end{enumerate}
\end{lemma}
\begin{proof}
    The required associativity and commutativity sequents follow from
    generalised identity, the two left-conjunction rules, and
    $\land\mathcal R$. Repeated use of these derivations, together with the
    admissibility of Cut, gives arbitrary changes of
    bracketing and order.  
    
    Clause (i) is obtained by repeated left-conjunction
    introduction around the chosen occurrence.  
    
    Clause (ii) follows from the
    conjunction projections in one direction and from $\land\mathcal R$ in the
    other. 
\end{proof}

The next theorem simultaneously prunes the derivation in \SGLDiamond and extract a \calc-derivation.

\begin{theorem}[Positive pruning and extraction]\label{thm:gl-positive-extraction}
    If $\SGLDiamond\vdash\Gamma\Rightarrow\Delta$, thene there is a formula occurrence $D$ in $\Delta$ such that both $\SGLDiamond\vdash\Gamma\Rightarrow D$ and $\calc\vdash\bigwedge\Gamma\implies D$.
\end{theorem}
\begin{proof}
    By induction on the height of the given $\SGLDiamond$-derivation, choosing the same succedent occurrence for the two claims.

    If the height is zero, the claims follows by the design of the systems and \Cref{lem:gl-conjunctive-contexts}(i).

    Suppose the last inference is $\land\mathcal L$:
        \begin{prooftree}
            \AxiomC{$A,B,\Gamma\Rightarrow\Delta$}
            \RightLabel{\scriptsize $\land\mathcal L$}
            \UnaryInfC{$A\land B,\Gamma\Rightarrow\Delta$}
        \end{prooftree}
    Then we apply the induction hypothesis to the premise and retain the selected
    occurrence $D$, using the design of \SGLDiamond and \Cref{lem:gl-conjunctive-contexts}(ii), by leveraging cut admissibility for \calc.

    Suppose the last inference is
    \begin{prooftree}
      \AxiomC{$\Gamma\Rightarrow\Delta,A$}
      \AxiomC{$\Gamma\Rightarrow\Delta,B$}
      \RightLabel{\scriptsize $\land\mathcal R$}
      \BinaryInfC{$\Gamma\Rightarrow\Delta,A\land B$}
    \end{prooftree}
    The induction hypotheses select occurrences $D_1$ in $\Delta,A$ and $D_2$
    in $\Delta,B$.  If $D_1$ is an occurrence in the displayed copy of
    $\Delta$, we retain $D_1$ and its two derivations; symmetrically, we do the same if
    $D_2$ lies in $\Delta$. If $D_1=A$ and $D_2=B$, then one application of $\land\mathcal R$ in $\SGLDiamond$ and one in \calc give the two required derivations with $D:=A\land B$.

    Finally, suppose the last inference is
    \begin{prooftree}
      \AxiomC{$A\Rightarrow\Sigma,\Diamond\Sigma,\Diamond A$}
      \RightLabel{\scriptsize $\Diamond$-L\"ob}
      \UnaryInfC{$\Omega,\Diamond A\Rightarrow\Diamond\Sigma,\Delta$}
    \end{prooftree}
    Apply the induction hypothesis to the premise and let $D'$ be the selected occurrence. By \Cref{lem:gl-positive-depth}, we have $D'\neq\Diamond A$ for the diagonal formula of the rule.

    If $D'=S$ is an occurrence in $\Sigma$, then $\SGLDiamond\vdash A\Rightarrow S$ and $\calc\vdash A\implies S$. By \Cref{lem:gl-positive-weakening}, we can extend the first derivation to $A\Rightarrow S,\Diamond S,\Diamond A$, and by applying $\Diamond$-L\"ob followed by left weakening we get $\SGLDiamond\vdash \Omega,\Diamond A\Rightarrow\Diamond S$. In \calc, Nec gives $\Diamond A\implies\Diamond S$, and \Cref{lem:gl-conjunctive-contexts}(i) adds the occurrences in $\Omega$. Thus we select the corresponding occurrence $D:=\Diamond S$ in $\Diamond\Sigma$.
    If $D'=\Diamond S$ is an occurrence in $\Diamond\Sigma$, the induction hypothesis gives
    $\SGLDiamond\vdash A\Rightarrow\Diamond S$ and $\calc\vdash A\implies\Diamond S$. As before, application of weakening yields $\SGLDiamond\vdash A\Rightarrow S,\Diamond S,\Diamond A$, from which we get $\SGLDiamond\vdash\Omega,\Diamond A\Rightarrow\Diamond S$ by applying $\Diamond$-L\"ob. In \calc, we apply 4Nec to $A\implies\Diamond S$ and get $\calc\vdash\Diamond A\implies\Diamond S$; \Cref{lem:gl-conjunctive-contexts}(i) adds $\Omega$.

    No other cases are left open, and we conclude the proof.
\end{proof}

\subsection{Syntactic proof of Dashkov's conservativity theorem}\label{subsec:gl-connection}

To achieve Dashkov's theorem, we only need an additional lemma.

\begin{lemma}[Embedding into $\gl$]\label{lem:calc-into-gl}
    If $\calc\vdash A\implies B$, then $\gl\vdash A\to B$.
\end{lemma}
\begin{proof}
    Straightforward induction on the height of the derivation in \calc.
\end{proof}

Collecting all the structural results we have proven, we can now state our central claim of conservativity over \gl.

\begin{theorem}[Dashkov's conservativity theorem]\label{thm:modular-dashkov}
    For all $A,B\in\lang$, the following equivalence holds:
    \[
        A\vdash_{\RC_1}B
    \quad\textit{iff}\quad
        \gl\vdash A\to B.
    \]
\end{theorem}
\begin{proof}
    If $A\vdash_{\RC_1}B$, then \Cref{theorem:CalcEquivToRCOne} gives
    $\calc\vdash A\implies B$, and \Cref{lem:calc-into-gl} yields
    $\gl\vdash A\to B$.

    Conversely, suppose $\gl\vdash A\to B$. By
    \Cref{cor:gl-positive-cut-free-completeness},
    $\SGLDiamond\vdash A\Rightarrow B$.  Applying
    \Cref{thm:gl-positive-extraction} to this singleton sequent yields
    $\calc\vdash A\implies B$ and  \Cref{theorem:CalcEquivToRCOne}
    then gives $A\vdash_{\RC_1}B$.
\end{proof}

\section{Uniform interpolation}
\label{sec:modular-uniform-interpolation}

Because of their minimalist nature, the calculi \calc and \qcalc makes their structural analysis quite transparent. In particular, uniform interpolation can be established in a rather natural way.  

In the quantified language, for instance, relation symbols,
individual constants, and free individual parameters can all be forgotten
by the same structural projection.  Bound variables are treated locally,
so no semantic construction is required.  The projection for $\RC_1$ reveals
then exactly the nullary special case.

Ordinary Craig interpolation for $\RC_1$ and $\QRC$ was already observed
in~\cite[Theorem~8.1]{BorgesJoosten:2026:SPLHA}: in the deductive
presentation, the non-logical symbols of a consequent are inherited from the
antecedent, so the consequent itself may serve as an interpolant.  

The result proven here is stronger in a different direction:  it is uniform in an
arbitrarily prescribed retained vocabulary and simultaneously eliminates
relation symbols, individual constants, and free individual parameters.
This places our construction within the broader proof-theoretic literature
on interpolation for modal and provability logics, including interpolation
for $\gl$ and $\mathsf{GLP}$ and general sequent-calculus methods for
uniform interpolation
\cite{Shamkanov:2011:Interpolation,Bilkova:2016:UniformInterpolation,
Iemhoff:2019:UniformInterpolation,TabatabaiJalali:2027:UniformInterpolation}.\footnote{
Recent constructive work also establishes uniform interpolation for modal
logics with an independent diamond by means of terminating, cut-free
sequent calculi~\cite{VanDerGiessenShillito:2026:UniformDiamond}.  In the
strictly positive setting considered here, the absence of negative polarity
makes the projection itself much simpler.}

\subsection{Left-uniform interpolant}

We write $\FV(A)$, $\Rel(A)$, $\Const(A)$, $\Tm(A)$, and $\Sub(A)$ for
the free variables, relation symbols, constants, terms, and subformulas
of $A$. The following simple lemma was essentially already observed in \cite{BorgesJoosten:2026:SPLHA}.

\begin{lemma}[Two invariants]
\label{lem:modular-invariants}
If
$\qcalc\vdash A\implies B$,
then
$\Rel(B)\subseteq\Rel(A)$
and
$\mdepth(B)\leq\mdepth(A)$.
\end{lemma}

\begin{proof}
Both claims follow by induction on the derivation, noticing that  substitution
changes terms but not relation symbols or modal depth, and that,  in the
$4\mathrm{Nec}$ case, the induction hypothesis gives the stronger
inequality $\mdepth(\Diamond B)\leq\mdepth(A)$.
\end{proof}

A simple corollary to this lemma is:

\begin{theorem}
    Any $\lang/\qlang$ formula is its own left-uniform interpolant in \calc/\qcalc respectively.
\end{theorem}

\begin{proof}
    Clearly, $B\implies B$; if $A \implies B$, then the common vocabulary is, in virtue of Lemma \ref{lem:modular-invariants}, exactly that of $B$. 
\end{proof}

Determining a right-uniform interpolant requires some more work.

\subsection{Full-vocabulary projection}
\label{subsec:full-vocabulary-projection}

A \emph{retained vocabulary} is a triple
$\Lambda=(\mathcal P,\mathcal C,X)$,
where $\mathcal P$ is a set of relation symbols, $\mathcal C$ a set of
individual constants, and $X$ a set of variables.  If $Y$ is a finite set
of variables currently bound above a displayed subformula, we define the
projection
$A\mathbin{\upharpoonright}_{\Lambda,Y}$
recursively as follows:  
\[
\begin{aligned}
\top\!\upharpoonright_{\Lambda,Y}
  &:= \top,\\
(B\land C)\!\upharpoonright_{\Lambda,Y}
  &:= (B\!\upharpoonright_{\Lambda,Y})
      \land (C\!\upharpoonright_{\Lambda,Y}),\\
(\Diamond B)\!\upharpoonright_{\Lambda,Y}
  &:= \Diamond(B\!\upharpoonright_{\Lambda,Y}),\\
(\forall x\,B)\!\upharpoonright_{\Lambda,Y}
  &:= \forall x\,
      (B\!\upharpoonright_{\Lambda,Y\cup\{x\}}).
\end{aligned}
\]
For an atom $S(t_1,\ldots,t_k)$, let
\[
  S(\vec t)\mathbin{\upharpoonright}_{\Lambda,Y}
  :=
  \begin{cases}
    S(\vec t),
      & \text{if }S\in\mathcal P,\text{ every constant in }\vec t
        \text{ lies in }\mathcal C,\\
      & \quad\text{and every variable in }\vec t
        \text{ lies in }X\cup Y;\\
    \top,
      & \text{otherwise.}
  \end{cases}
\]
Finally, write
$A\mathbin{\upharpoonright}_{\Lambda} :=A\mathbin{\upharpoonright}_{\Lambda,\varnothing}$.\footnote{
The definition is invariant under $\alpha$-conversion: renaming a binder
renames both its bound occurrences and the corresponding member of $Y$.}

A formula $B$ will be called a $(\Lambda,Y)$-formula if
\[
  \Rel(B)\subseteq\mathcal P,
  \qquad
  \Const(B)\subseteq\mathcal C,
  \qquad
  \FV(B)\subseteq X\cup Y.
\]
In particular,
$A\mathbin{\upharpoonright}_{\Lambda,Y}$ is a
$(\Lambda,Y)$-formula.

We first record two elementary syntactic facts about the projection.

\begin{lemma}[Projection and substitution]
\label{lem:projection-substitution}
Let $x\notin Y$, and let $t$ be free for $x$ in $A$.  After
standardising bound variables apart from $Y\cup\{x\}\cup\FV(t)$,
\[
  \Bigl(
    (A\mathbin{\upharpoonright}_{\Lambda,Y\cup\{x\}})
    \subst{x}{t}
  \Bigr)
  \mathbin{\upharpoonright}_{\Lambda,Y}
  \quad\text{is $\alpha$-equivalent to}\quad
  (A\subst{x}{t})\mathbin{\upharpoonright}_{\Lambda,Y}.
\]
Moreover, if $x\notin\FV(A)$, then 
  $A\mathbin{\upharpoonright}_{\Lambda,Y\cup\{x\}}$
  is $\alpha$-equivalent to 
  $A\mathbin{\upharpoonright}_{\Lambda,Y}$.
\end{lemma}

\begin{proof}
For each formula $A$, we prove both statements simultaneously for every
finite $Y$ satisfying the displayed side conditions, by induction on $A$.

Only the atomic retention test requires a comment.  An occurrence of $x$ is permitted on the left
because $x\in Y\cup\{x\}$; after substitution, the resulting atom is
retained by the outer projection precisely when every variable occurring
in $t$ is allowed by $X\cup Y$ and every constant occurring in $t$ is
retained.  This is exactly the retention test applied directly to
$A\subst{x}{t}$.  

The quantifier case uses the standing $\alpha$-conversion convention to keep the displayed binder away from
$t$ and the finite parameter set under consideration.
\end{proof}

\begin{lemma}
\label{lem:projection-derivable}
For every $A\in\qlang$ and every finite $Y$,
$\qcalc\vdash A\implies A\mathbin{\upharpoonright}_{\Lambda,Y}$.
\end{lemma}

\begin{proof}
We perform an induction on $A$, with the strengthened statement quantified over every
finite $Y$.  

For an atom, use Initial when it is retained and
$\top\mathcal R$ when it is erased.  The cases for $\top$, conjunction,
and $\Diamond$ are reconstructed with the corresponding rules.

Let then $A=\forall x\,C$.  By inductive hypothesis,
$\qcalc\vdash C\implies C\mathbin{\upharpoonright}_{\Lambda,Y\cup\{x\}}$.
Applying $\forall\mathcal L$ with the term $x$ gives
$\qcalc\vdash \forall x\,C \implies C\mathbin{\upharpoonright}_{\Lambda,Y\cup\{x\}}$.
Since $x\notin\FV(\forall x\,C)$, $\forall\mathcal R$ yields
  \mbox{$\qcalc\vdash
  \forall x\,C
  \implies
  \forall x\,
  (C\mathbin{\upharpoonright}_{\Lambda,Y\cup\{x\}})$},
which is the required projection.
\end{proof}

The next theorem establishes the universal property of the construction.

\begin{theorem}[Full-vocabulary uniform projection]
\label{thm:full-vocabulary-projection}
    Let $Y$ be finite and let $B$ be a $(\Lambda,Y)$-formula.  
    
    If
    $\qcalc\vdash A\implies B$,
    then
    $\qcalc\vdash A\mathbin{\upharpoonright}_{\Lambda,Y}\implies B$.
\end{theorem}

\begin{proof}
By cut elimination (\Cref{thm:qcalc-cut}), we can perform an induction on a cut-free derivation, with the induction statement
strengthened simultaneously over every finite set $Y$: whenever the
end-sequent is $A\implies B$ and $B$ is a $(\Lambda,Y)$-formula,
derive $A\mathbin{\upharpoonright}_{\Lambda,Y}\implies B$.\footnote{This strengthening is needed only because a right-universal step passes
from $Y$ to $Y\cup\{x\}$.}

Initial is immediate because an atomic $(\Lambda,Y)$-consequent passes
the retention test, and $\top\mathcal R$ is immediate as well.  

The conjunction and modal cases are obtained by applying the induction
hypothesis to the premise (or premises) and reconstructing the last rule.

Suppose the last inference is $\forall\mathcal R$:
\begin{prooftree}
  \AxiomC{$A\implies C$}
  \RightLabel{\scriptsize $\forall\mathcal R$}
  \UnaryInfC{$A\implies\forall x\,C$}
\end{prooftree}
with $x\notin\FV(A)$.  Since $\forall x\,C$ is a
$(\Lambda,Y)$-formula, $C$ is a
$(\Lambda,Y\cup\{x\})$-formula.  The induction hypothesis gives
$A\mathbin{\upharpoonright}_{\Lambda,Y\cup\{x\}} \implies C$.
By Lemma~\ref{lem:projection-substitution} and
$x\notin\FV(A)$, its antecedent is $\alpha$-equivalent to
$A\mathbin{\upharpoonright}_{\Lambda,Y}$, and we can reapply $\forall\mathcal R$.

Suppose the last inference is $\forall\mathcal L$:
\begin{prooftree}
  \AxiomC{$C\subst{x}{t}\implies B$}
  \RightLabel{\scriptsize $\forall\mathcal L$}
  \UnaryInfC{$\forall x\,C\implies B$}
\end{prooftree}
where $t$ is free for $x$ in $C$.  Standardising the displayed binder so
that $x\notin Y$ and the substitution is capture-avoiding,  the inductive hypothesis gives
$(C\subst{x}{t}) \mathbin{\upharpoonright}_{\Lambda,Y} \implies B$.
By Lemma~\ref{lem:projection-derivable} applied after substitution, and
then Lemma~\ref{lem:projection-substitution}, we have
  $(C\mathbin{\upharpoonright}_{\Lambda,Y\cup\{x\}})
  \subst{x}{t}
  \implies
  (C\subst{x}{t})
  \mathbin{\upharpoonright}_{\Lambda,Y}$.
Cut admissibility therefore
gives
$(C\mathbin{\upharpoonright}_{\Lambda,Y\cup\{x\}}) \subst{x}{t} \implies B$.
One application of $\forall\mathcal L$ yields
$\forall x\, (C\mathbin{\upharpoonright}_{\Lambda,Y\cup\{x\}}) \implies B$,
which is exactly
$(\forall x\,C)\mathbin{\upharpoonright}_{\Lambda,Y} \implies B$.
This completes the proof.
\end{proof}

\subsection{Uniform and Craig interpolation}
\label{subsec:modular-interpolation-consequences}

We are ready to see how the previous lemmas relate to the interpolation analysis we aim to.  
The preceding theorem says that the projection is the strongest consequence of $A$ expressible in the retained vocabulary, taking $Y:=\varnothing$.

\begin{theorem}[Uniform post-interpolation for \QRC]
\label{thm:qrc-uniform-interpolation}
Let $A\in\qlang$ and let
$\Lambda=(\mathcal P,\mathcal C,X)$.  Let us define.
$U_{\Lambda}(A) :=A\mathbin{\upharpoonright}_{\Lambda}$.
Then the following hold:
\begin{enumerate}
  \item $\qcalc\vdash A\implies U_{\Lambda}(A)$;
  \item $U_{\Lambda}(A)$ is a $\Lambda$-formula;
  \item for every $\Lambda$-formula $B$,
    $\qcalc\vdash A\implies B$
    iff
    $\qcalc\vdash U_{\Lambda}(A)\implies B$.
\end{enumerate}
\end{theorem}

\begin{proof}
The first clause is Lemma~\ref{lem:projection-derivable}; the second is
immediate from the definition.  The forward direction of the third clause
is Theorem~\ref{thm:full-vocabulary-projection}.  For the converse, compose
$A\implies U_{\Lambda}(A)$ with
$U_{\Lambda}(A)\implies B$ using cut admissibility.
\end{proof}

By Theorem~\ref{theorem:QRC1andQlangEquivalent}, the same statement
holds considering the standard axiomatisation of \QRC.

\begin{corollary}[Craig interpolation for \QRC]
\label{cor:qrc-craig-interpolation}
If
$A\vdash_{\QRC}B$,
then there is a formula $I$ such that
\begin{enumerate}
    \item $A\vdash_{\QRC}I$ and $I\vdash_{\QRC}B$;
    \item $\Rel(I)\subseteq\Rel(A)\cap\Rel(B)$,
            $\Const(I)\subseteq\Const(A)\cap\Const(B)$, and
            $\FV(I)\subseteq\FV(A)\cap\FV(B)$.
\end{enumerate}
\end{corollary}

\begin{proof}
Take
$\Lambda=(\Rel(B),\Const(B),\FV(B))$
and put $I:=U_{\Lambda}(A)$.  The two derivabilities follow from
Theorem~\ref{thm:qrc-uniform-interpolation}.  Since projection never adds a
relation symbol, constant, or free variable, the three displayed
inclusions follow immediately.
\end{proof}

Since every occurrence in the strictly positive language has positive
polarity, the relation-symbol part of the preceding result gives also the
Lyndon property of interpolants.

\subsection{The nullary specialisation}
\label{subsec:rc-uniform-specialisation}

For $\RC_1$, we can identify each propositional variable with a nullary relation
symbol.  There are no individual terms, so the constant and free-parameter
components of a retained vocabulary disappear.  For a set $\Gamma$ of
propositional variables, the projection discussed above reduces to
\[
  U_{\Gamma}(\top)=\top,
  \qquad\text{and}\qquad
  U_{\Gamma}(p)=
  \begin{cases}
    p,&p\in\Gamma,\\
    \top,&p\notin\Gamma,
  \end{cases}
\]
with $U_{\Gamma}$ commuting with $\land$ and $\Diamond$.

Interpolation for the propositional reflection calculus easily follows then.

\begin{corollary}[Uniform post-interpolation for $\RC_1$]
\label{cor:rc-uniform-interpolation}
For every $A\in\lang$ and every set $\Gamma$ of propositional variables,
$U_{\Gamma}(A)$ satisfies the following:
\begin{enumerate}
    \item $A\vdash_{\RC_1}U_{\Gamma}(A)$;
    \item $\operatorname{Var}(U_{\Gamma}(A))\subseteq\Gamma$; and
    \item for every $B$ with $\operatorname{Var}(B)\subseteq\Gamma$,
 $A\vdash_{\RC_1}B$ iff $U_{\Gamma}(A)\vdash_{\RC_1}B$.
\end{enumerate}
That is: the uniform interpolation theorem for $\RC_1$ is exactly the
nullary instance of Theorem~\ref{thm:qrc-uniform-interpolation}.
\end{corollary}

\begin{proof}
Specialise Theorem~\ref{thm:qrc-uniform-interpolation} to nullary relation
symbols and use deductive conservativity of \qcalc over \calc.
\end{proof}

To conclude, we remark that the projection we defined is size-preserving at the level of logical connectives and
can only shorten atoms.  With constant-time membership tests for the
retained vocabulary, it is computed in $\mathcal O(|A|)$ time and satisfies
\[
  |U_{\Lambda}(A)|\leq |A|,
  \qquad
  \mdepth(U_{\Lambda}(A))=\mdepth(A),
  \qquad
  \qdepth(U_{\Lambda}(A))=\qdepth(A).
\]
Thus the quantified strengthening does not sacrifice the elementary,
constructive character of the propositional projection.

\section{Conclusions and future work}\label{sec:conclusion}

We have presented novel cut-free proof systems for the strictly positive fragment of Gödel-Löb's provability logic $\gl$~\cite{Boolos:1993:LogicOfProvability,sep-logic-provability} and its universal quantified extension. Throughout this work, we denote the propositional strictly positive fragment by $\RC_1$, and its quantified lifting by \QRC.

The three sequent calculi (\calc, \calcFourNec, and $\calcFourNec_*$) proposed here are designed as follows: 
\begin{enumerate}
    \item \calc~includes an analytic modal rule for necessitation and an analytic modal rule for the transitivity of diamonds;
    \item \calcFourNec~preserves the modal rule for necessitation but replaces the analytic rule for transitivity with a pair of rules -- these handle transitivity of diamonds analytically on the left-hand side of sequents, but non-analytically on the right;
    \item $\calcFourNec_*$ retains the analytic rule for necessitation and the non-analytic rule for transitivity of diamonds on the right-hand side, whilst removing the companion rule for diamonds on the left.
\end{enumerate} 

Their first-order version (\qcalc, \qcalcFourNec, and $\qcalcFourNec_*$) are obtained by adding a natural module for $\forall$-rules on sequents.

Even for the non-analytic calculi, as main theorems, we have proven that all our sequent systems admit an effective cut elimination procedure. Such a structural result, let us develop purely syntactic and constructive proofs of known though non-trivial results about $\RC_1$ and \QRC, namely: conservativity over \gl (Dashkov's theorem) and (modular) uniform interpolation.

\medskip

A non-nested cut-free calculus for $\RC_\Lambda$ (where $\Diamond$-operators are indexed over ordinals in $\Lambda$) or the impossibility thereof remains an open problem.

We may consider the following rules
   \begin{prooftree}
    \AxiomC{$A \implies B$}\RightLabel{\scriptsize MonoNec \ \ \ \ \ where $\beta \leq \alpha$}\UnaryInfC{$ \langle \alpha \rangle A \implies \langle \beta \rangle B$}
    \end{prooftree}
and
    \begin{prooftree}
        \AxiomC{$ A \implies \langle \alpha \rangle B$}
    \RightLabel{\scriptsize OrdinalNec}
    \UnaryInfC{$ \langle \beta \rangle A \implies \langle \min \{\alpha,\beta\} \rangle B$}
    \end{prooftree}

By adding these rules to the diamond-free fragment of \calc we easily get the following.

\begin{theorem}
    The resulting calculus generates a fragment of $\RC_\Lambda$ and has cut-admissibility.
\end{theorem}

The proof is very similar to that of \Cref{theorem:CutAdmissibilityQRCOne,thm:calc-cut}.
However, it is also easy to see that the calculus cannot prove Negative introspection, that is, 
$\langle \alpha \rangle A \wedge \langle \beta \rangle B \implies \langle \alpha \rangle (A \wedge \langle \beta \rangle B)$ for $\beta < \alpha < \Lambda$.

How can this calculus be enlarged while retaining cut-admissibility? Also, does the calculus correspond to the strictly positive fragment of some natural polymodal provability logic? 

Related questions are if $\QRC$ corresponds in the same sense to the quantified version of $\gl$. Polymodal versions of $\QRC$ would be very natural to study and, potentially they could be useful in the study of ordinal analysis through iterated reflection (see e.g.~\cite{Beklemishev:2005:Survey}). Providing well-behaved relational semantics for them have proven to be a non-trivial generalisation of Kripke semantics for $\QRC$. Good calculi could thus provide well-behaved model extractions for unprovable sequents.

\paragraph{Acknowledgments} 
The first author has received financial support through the project project SERICS – Security and Rights in the CyberSpace PE0000014 within PNRR, M4C2 I.1.3, funded by the European Union - NextGenerationEU (MUR Code: 2022CY2J5S). We thank the Istituto Nazionale di Alta Matematica/INdAM group GNSAGA for institutional support.

The second author has received financial support through the projects PID2023-151396OB-I00 and PID2023-149556NB-I00 of the Spanish Ministry of Science and Innovation, through 2021 SGR 00348 of the \textit{Generalitat de Catalunya} and, through the Severo Ochoa and María de Maeztu Program for Centers and Units of Excellence in R\&D (CEX2020-001084-M), the Spanish State Research Agency. We thank CERCA Programme/Generalitat de Catalunya for institutional support.

\printbibliography

\end{document}